\documentclass{agtart_a} 
\pdfoutput=1

\usepackage{euscript}           % Nice script font.
\usepackage{enumerate,calc}
\usepackage[matrix,arrow,curve,frame]{xy}    % XY-pic diagram pac

\title{Postnikov extensions of ring spectra}

\author{Daniel Dugger}
\givenname{Daniel}
\surname{Dugger}
\address{Department of Mathematics\\ 
University of Oregon\\\newline
Eugene, OR 97403\\USA}
\email{ddugger@math.uoregon.edu}
\urladdr{}

\author{Brooke Shipley}
\givenname{Brooke}
\surname{Shipley}
\address{Department of Mathematics\\
University of Illinois at Chicago\\\newline
Chicago, IL 60607\\USA}
\email{bshipley@math.uic.edu}
\urladdr{}

\volumenumber{6}
\issuenumber{}
\publicationyear{2006}
\papernumber{63}
\startpage{1785}
\endpage{1829}

\doi{}
\MR{}
\Zbl{}


\received{26 July 2006}
\revised{}
\accepted{22 August 2006}
\published{1 November 2006}
\publishedonline{1 November 2006}
\proposed{}
\seconded{}
\corresponding{}
\editor{}
\version{}

\arxivreference{math.AT/0604260}

\let\xysavmatrix\xymatrix
\def\xymatrix{\disablesubscriptcorrection\xysavmatrix}
\AtBeginDocument{\let\bar\wbar\let\tilde\wtilde\let\hat\what}
\newcommand{\dfn}{\textit} % Make defined words bold.
\newcommand{\mdfn}[1]{\dfn{#1\/}} % Even make math bold

\newcommand{\Alg}{\cA lg}
\makeop{Id}                           % The identity functor
\newcommand{\id}{\Id}                            % The identity functor
\makeop{PB}
\newcommand{\we}{\smash{\rlap{\kern 6pt\raise 4pt\hbox{\footnotesize $\sim$}}}\longrightarrow}
\newcommand{\bwe}{\smash{\rlap{\kern 8.5pt\raise 4pt\hbox{\footnotesize $\sim$}}}\longleftarrow}
\newcommand{\trfib}{\smash{\rlap{\kern 7pt\raise 4pt\hbox{\footnotesize $\sim$}}}\longfib}
\newcommand{\trcof}{\smash{\rlap{\kern 5.5pt\raise 4pt\hbox{\footnotesize $\sim$}}}\rightarrowtail}

\xymatrixcolsep{1.9pc}                          % Adjust size of diagrams.
\xymatrixrowsep{1.9pc}
\newdir{ >}{{}*!/-8pt/\dir{>}}                  % Make better tailed arrows

\def\longfib{\DOTSB\relbar\joinrel\twoheadrightarrow}

\makeatletter
\def\cnewtheorem#1[#2]#3{\newtheorem{#1}{#3}[section]
\expandafter\let\csname c@#1\endcsname\c@subsection}
\makeatother

\swapnumbers
\cnewtheorem{thm}[subsection]{Theorem}
\cnewtheorem{prop}[subsection]{Proposition}
\cnewtheorem{claim}[subsection]{Claim}
\cnewtheorem{cor}[subsection]{Corollary}
\cnewtheorem{lemma}[subsection]{Lemma}

\theoremstyle{definition}  % Bold headings and Roman body text.
\cnewtheorem{defn}[subsection]{Definition}
\cnewtheorem{example}[subsection]{Example}
\cnewtheorem{exercise}[subsection]{Exercise}
\cnewtheorem{note}[subsection]{Note}
\cnewtheorem{warning}[subsection]{Warning}
\cnewtheorem{remark}[subsection]{Remark}
\cnewtheorem{discussion}[subsection]{Discussion}

\numberwithin{equation}{subsection}
\newcommand{\Smash}             {\wedge}
\newcommand{\sm}             {\wedge}

\newcommand{\Wedge}             {\vee}

\newcommand{\tens}              {\otimes}               %tensor
\newcommand{\iso}               {\cong}
\newcommand{\cat}{\EuScript}    % Use \EuScript to name a category.
\newcommand{\cA}{{\cat A}}      % Only seems to work for caps, and only gets
\newcommand{\cB}{{\cat B}}      % first letter.
\newcommand{\cC}{{\cat C}}
\newcommand{\CC}{F_C}
\newcommand{\cD}{{\cat D}}
\newcommand{\cE}{{\cat E}}

\newcommand{\ccH}{{\cat Hom}}

\newcommand{\cK}{{\cat K}}

\newcommand{\cM}{{\cat M}}

\newcommand{\Ho}{\text{Ho}\,}
\newcommand{\ho}{\text{Ho}\,}

\newcommand{\field}[1]  {\mathbb #1} % Use blackboard bold for these sets

\DeclareMathOperator*{\colim}{colim}

\DeclareMathOperator{\Hom}{Hom}

\DeclareMathOperator{\Ext}{Ext}
\DeclareMathOperator{\Tor}{Tor}

\DeclareMathOperator{\Map}{Map}

\DeclareMathOperator{\Der}{Der}

\DeclareMathOperator{\THH}{THH}
\DeclareMathOperator{\Aut}{Aut}
\DeclareMathOperator{\hAut}{hAut}

\DeclareMathOperator{\Cell}{--Cell}
\DeclareMathOperator{\bimod}{--bimod}
\DeclareMathOperator{\lmod}{--Mod}

\DeclareMathOperator{\Cbimod}{C--bimod}

\newcommand{\ra}{\rightarrow}                   % right arrow
\newcommand{\lra}{\longrightarrow}              % long right arrow
\newcommand{\la}{\leftarrow}                    % left arrow
\newcommand{\from}{\leftarrow}                    % left arrow
\newcommand{\lla}{\longleftarrow}               % long left arrow
\newcommand{\llra}[1]{\smash{\stackrel{#1}{\lra}}}      % labeled long right
\newcommand{\llla}[1]{\smash{\stackrel{#1}{\lla}}}      % labeled long right
\newcommand{\weak}{homotopy }                   % weak equivalence
\newcommand{\cof}{\rightarrowtail}              % cofibration
\newcommand{\fib}{\twoheadrightarrow}           % fibration

\newcommand{\inc}{\hookrightarrow}              % inclusion
\newcommand{\binc}{\hookleftarrow}              % inclusion
\newcommand{\blank}{-}                          % A hyphen, as in
\newcommand{\mC}{\underline{\cC}}

\newcommand{\ovcat}{\downarrow}

\newcommand{\bd}{\partial}

\newcommand{\adjoint}{\rightleftarrows}

\newcommand{\bdd}[1]{\partial\Delta^{#1}}
\newcommand{\del}[1]{\Delta^{#1}}

\newcommand{\he}{\simeq}

\newcommand{\PP}{{\field P}}   %Ordinary Postnikov functor

\newcommand{\NonU}{{\cat NonU}}

\renewcommand{\theequation}{\thesubsection}

\newenvironment{myequation}
  {\addtocounter{subsection}{1}\begin{equation}}
  {\end{equation}}

\newcommand{\lMod}{\lmod}

\newcommand{\mc}{\co}

\DeclareMathOperator{\Ev}{Ev}

\newcommand{\cofib}{\cof}
\newcommand{\Ring}{\cat{R}ing}

\newcommand{\barQ}{\underline{Q}}
\newcommand{\barZ}{\underline{Z}}

\newcommand{\bfib}{\twoheadleftarrow}           % fibration
\newcommand{\bcof}{\leftarrowtail}              % cofibration

\newcommand{\tT}{\widetilde{T}_C}

\newcommand{\fC}{fC}

\newcommand{\fP}{fPC}
\newcommand{\coprd}{\amalg}

\begin{document}

\begin{asciiabstract}
We give a functorial construction of k-invariants for ring spectra
and use these to classify extensions in the Postnikov tower of a ring
spectrum.  
\end{asciiabstract}

\begin{htmlabstract}
We give a functorial construction of k&ndash;invariants for ring spectra
and use these to classify extensions in the Postnikov tower of a ring
spectrum.  
\end{htmlabstract}

\begin{abstract}
We give a functorial construction of $k$--invariants for ring spectra
and use these to classify extensions in the Postnikov tower of a ring
spectrum.  
\end{abstract}

\maketitle

\section{Introduction}

This paper concerns $k$--invariants for ring spectra and
their role in classifying Postnikov extensions.
Recall that a connective ring spectrum $R$ has a Postnikov tower
\[ \cdots \ra P_2R \ra P_1R \ra P_0R \ra *
\]
in the homotopy category of ring spectra.  The levels come equipped
with compatible maps $R\ra P_n R$, and the $n$--th level is
characterized by having $\pi_i(P_nR)=0$ for $i>n$, together with the
fact that $\pi_i(R) \ra \pi_i(P_nR)$ is an isomorphism for $i\leq n$.
In this paper we produce $k$--invariants for the levels of this tower
and explain their role in the following problem: if one only knows
$P_{n-1}R$ together with $\pi_n(R)$ as a $\pi_0(R)$--bimodule, what are
the possibilities for $P_nR$?  \fullref{cor:main} shows
in what sense the possibilities are classified by $k$--invariants.

\subsection[Classical k-invariants]{Classical $k$--invariants}
To explain our results further, it's useful to briefly recall the
situation for ordinary topological spaces.  If $X$ is a space, let
$P_nX$ be the $n$--th Postnikov section of $X$.  The
$k$--invariant is a map $P_{n-1} X \ra K(\pi_{n}X,n+1)$ and the
homotopy fiber of this map is weakly equivalent to $P_{n}X$.  So
$P_nX$ can be recovered from the $k$--invariant, and in fact the
$k$--invariant only depends on $P_nX$.
One
is tempted to say that the possibilities for $P_nX$ are {\it
classified by\/} the possible $k$--invariants, but this is where some
care is needed.

To clarify the situation, it's useful to set $C=P_{n-1}X$ and
$M=\pi_{n}X$.  By a Postnikov $n$--extension of $C$ (of type $M$) we
mean a space $Y$ together with a map $Y\ra C$ such that $\pi_i(Y)\ra
\pi_i(C)$ is an isomorphism for $i\leq n-1$, $\pi_n Y\iso M$, and
$\pi_{i}Y=0$ for $i>n$.  Note that the isomorphism $\pi_nY\iso M$ is
not part of the data.  The Postnikov $n$--extensions form a category,
in which a map from $Y\ra C$ to $Y'\ra C$ is a weak equivalence $Y\ra
Y'$ making the evident triangle commute.  We'll denote this category
$\cM(C,M,n)$.

For convenience, suppose $C$ is simply connected (so we
don't have to worry about the $\pi_1C$ actions).  One is tempted to
claim that the connected components of $\cM(C,M,n)$ are in bijective
correspondence with the set of homotopy classes $[C,K(M,n+1)]$.
Unfortunately, this isn't quite the case.  Note that the group $\Aut
M$ of abelian group automorphisms acts on $K(M,n+1)$ and hence on
$[C,K(M,n+1)]$.  If a certain $k$--invariant $C\ra K(M,n+1)$ is
``twisted'' by an automorphism of $M$, it gives rise to a weakly
equivalent extension of $C$.  The correct statement, it then turns
out, is that if $C$ is simply connected there is a bijection
\[ \pi_0\cM(C,M,n) \iso [C,K(M,n+1)]/\Aut(M).
\]
This statement is best proven by upgrading it to a statement about the
homotopy type of $\cM(C,M,n)$ (where by the homotopy type of a
category we always mean the homotopy type of
its nerve).  One can prove that there is a homotopy fiber sequence
\[ \Map(C,K(M,n+1)) \ra \cM(C,M,n) \ra B\Aut(M),
\]
and the resulting long exact sequence of homotopy groups gives the
identification of $\pi_0\cM(C,M,n)$ cited above.

For a proof of this homotopy fiber sequence (together with a version
when $C$ is not simply connected) we refer the reader to
Blanc, Dwyer and Goerss \cite[Sections 2,3]{BDG}.  An important part of the proof is having a
simple, functorial construction of the $k$--invariant, and we now
describe this.  If $p\co Y\ra C$ is a Postnikov $n$--extension of
type $M$, let $D$ denote the homotopy cofiber of $p$.  One can prove
by a Blakers--Massey type result that $\pi_i(D)=0$ for $i\leq n$,
whereas $\pi_{n+1}D \iso \pi_n Y\iso M$.  Then $P_{n+1}D$ is an
Eilenberg--Mac Lane space $K(M,n+1)$, and our $k$--invariant for $Y$ is
the composite
\[ C \ra D \ra P_{n+1}D.
\]
Note that in some sense this is not really a $k$--invariant, as
one does not have a specified weak equivalence $P_{n+1}D\he
K(M,n+1)$.  One can prove that different weak equivalences differ by
an element of $\Aut(M)$, and this shows that one has a well-defined
element of the orbit space $[C,K(M,n+1)]/\Aut(M)$.

\subsection{Results for ring spectra}
Now we jump into the category of ring spectra, and state our main
results.  In this paper we always work in the category of symmetric
spectra from \cite{hss}.  So ``ring spectrum'' means ``symmetric ring
spectrum''.

Fix $n\geq 1$.  Let $C$ be a connective ring spectrum such
that $P_{n-1}C\he C$, and let $M$ be a $\pi_0C$--bimodule.  By a
\mdfn{Postnikov extension of $C$ of type $(M,n)$} one means a 
ring map $Y\ra C$ such that
\begin{enumerate}[(i)]
\item $\pi_iY=0$ for $i>n$,
\item $\pi_i Y\ra \pi_i C$ is an isomorphism for $i\leq n-1$,
\item $\pi_n Y\iso M$ as $\pi_0Y$--bimodules (where $M$ becomes a
$\pi_0 Y$--bimodule via the isomorphism $\pi_0 Y\ra \pi_0C$).
\end{enumerate}
A map of Postnikov extensions from $Y\ra C$ to $Y'\ra C$ is a weak
equivalence of ring spectra $Y\ra Y'$ making the triangle commute.
Denote the resulting category by $\cM(C{+}(M,n))$.  We call this the
\dfn{moduli space\/} (or \dfn{moduli  category\/}) of Postnikov extensions
of $C$ of type $(M,n)$.  Note that we are not assuming that $C$ is
fibrant here.

We next identify the analogues of Eilenberg--Mac Lane spaces.
Given a $C$--bimodule $W$, one can construct a ring spectrum $C\Wedge
W$ whose underlying spectrum is the wedge and where the multiplication
comes from the bimodule structure on $W$---so $W$ squares to
zero under this product.  We call $C\Wedge W$ the \dfn{trivial square
zero extension} of $C$ by $W$.

Given our $\pi_0(C)$--bimodule $M$, there is a $C$--bimodule $HM$ for
which $\pi_0(HM)\iso M$ and $\pi_i(HM)=0$ for $i\neq 0$.  In fact, all
such bimodules are weakly equivalent (in the category of
$C$--bimodules).  One gets resulting bimodules $\Sigma^i(HM)$ for all
$i$, and therefore ring spectra $C\Wedge \Sigma^i (HM)$.  Throughout
the paper we will abuse notation and simplify $HM$ to just
$M$---thus, we will write $C\Wedge \Sigma^i M$ for $C\Wedge \Sigma^i
HM$.

Here is our main theorem:
\begin{thm}
\label{thm:main}
Given $C$ and $M$ as above, there is a homotopy fiber sequence
\[ \underline{\Ring_{\smash{/C}}}(C,C\Wedge \Sigma^{n+1} M)
\ra \cM(C{+}(M,n)) \ra B\Aut(M)
\]
where $\underline{\Ring_{\smash{/C}}}(X,Y)$ denotes the homotopy mapping
space from $X$ to $Y$ in the category of ring spectra over $C$ and 
$\Aut(M)$ is the group of $\pi_0(C)$--bimodule automorphisms of
$M$.
\end{thm}

\begin{cor}
\label{cor:main}
There is a bijection of sets
\[ \pi_0\cM(C{+}(M,n)) \iso [\Ho(\Ring_{/C})(C,C\Wedge \Sigma^{n+1}
M)]/\Aut(M).
\]
\end{cor}

In the context of these results, the main difference between ring
spectra and ordinary topological spaces is that there are no absolute
cohomology theories for ring spectra.  When dealing with ring spectra,
one always deals with {\it relative\/} cohomology theories.  (For a
nice explanation of this phenomenon in the commutative case, see the
introduction to Basterra and Mandell \cite{BM}.)  Thus, in the above results one is forced
to always work over $C$: the analogue of the Eilenberg--Mac Lane space is
the ring spectrum $C\Wedge \Sigma^{n+1}M$, and the mapping spaces must
be computed in the category of ring spectra {\it over $C$\/}.  Aside
from these differences, the statements for ring spectra and
topological spaces are quite similar.

The set of homotopy classes appearing in \fullref{cor:main} can
be identified with a topological Hochschild cohomology group.  In this
way one sees that $\THH^*$ is the natural receptor for $k$--invariants of
associative ring spectra.  This perspective also simplifies
calculations, since topological Hochschild cohomology involves a
homotopy category of bimodules instead of a homotopy category of
rings.  All of this is discussed in \fullref{sec:thh}, and we have
the following restatement of \fullref{cor:main}.

\begin{prop}\label{prop:main}
Let $C$ and $M$ be as above. 
Assuming that $C$ is cofibrant as an underlying spectrum, one has a bijection 
\[\pi_0\cM(C{+}(M,n)) \iso \THH^{n+2}(C,M)/\Aut(M).\]
\end{prop}

Finally, we remark that
the above results can actually be extended, so that they apply not
just to ring spectra but to algebras over a given connective,
commutative ring spectrum $R$.  This is the form in which we will
actually prove them (see \fullref{thm:classify} and
\fullref{prop-8.8}).  Moreover, all the
results of the paper apply equally well to the category of
differential graded algebras over a commutative ground ring $k$.  Our
proofs all adapt essentially verbatim, or else one can use that the
homotopy theory of dgas over $k$ is equivalent to that of algebras
over the Eilenberg--Mac Lane ring spectrum $Hk$ (this is proven
by the second author in \cite{S}).    

\begin{remark}
Both \fullref{thm:main} and \fullref{cor:main} can also be
rewritten in terms of ring spectra over $P_0C$ (the zero-th Postnikov
stage of $C$) rather than ring spectra over $C$.  
For this, 
see \fullref{prop-PC}.
\end{remark}

\subsection{Some background}
\fullref{cor:main} was needed in our paper \cite{tpwe}, and we
at first believed this result to be obvious.  Our attempts to give a
careful proof, however, always seemed to fail.  A construction of
$k$--invariants for ring spectra had already been given in Lazarev \cite{L},
but that construction does not seems well-suited for the above
classification questions.  A construction of $k$--invariants for {\it
commutative\/} ring spectra appeared in Basterra \cite{B}; while this
construction also did not meet all of our needs, many of the
techniques of \cite{B} are used in our \fullref{sec:nonu}.
\vspace{-3pt}

Eventually we discovered Blanc, Dwyer and Goerss \cite{BDG}, which applied the Dwyer--Kan
moduli space technology to the classification of Postnikov extensions
in a related context.  It will be clear to the reader that the basic
methods in the present paper are heavily influenced by \cite{BDG}
(there is one main difference, discussed in 
\fullref{re:bdg}).  However, in order
to carry out the program of \cite{BDG} we have had to straighten out many
points about ring spectra along the way.  One of the main things the
reader will find here is a careful, functorial construction of
$k$--invariants for ring spectra.  We also provide a careful proof of a
Blakers--Massey theorem for ring spectra in \fullref{se:pushout-proof}.  We thank a
helpful referee and Mike Mandell for suggestions which improved our
presentation of this result.
\vspace{-3pt}

\subsection{Notation and terminology}
If $\cC$ is a category then we write $\cC(X,Y)$ for $\Hom_\cC(X,Y)$.
If $\cC$ is a simplicial model category we write $\Map_{\cC}(X,Y)$ for
the simplicial function complex.  If $\cC$ is a model category then
$\underline{\cC}(X,Y)$ denotes a homotopy function complex from $X$ to
$Y$.  The phrase ``homotopy function complex'' indicates a construction
which has the correct homotopy type even if $X$ is not cofibrant and
$Y$ is not fibrant.  To fix a particular construction, we use the
hammock localization of Dwyer and Kan \cite{function}.
\vspace{-3pt}

If $R$ is a commutative ring spectrum then there are model category
structures on $R\lMod$ and $R$--$\Alg$ provided by Schwede and the second author \cite{SS1};
in each case the fibrations and weak equivalences are
determined by the forgetful functor to symmetric spectra.  We use
these model categories throughout the paper.
\vspace{-3pt}

To every category $\cC$ one can associate a simplicial set $|\cC|$ by
taking its nerve.  In this paper we often abbreviate $|\cC|$ to just
$\cC$, letting the application of the nerve be clear from context.
\vspace{-3pt}

Finally, if $X$ is a spectrum then $\pi_* X$ always refers to the derived
homotopy groups (that is, homotopy groups of a fibrant replacement).

\subsubsection{Acknowledgments}
%\thanks{Second author partially supported by NSF Grant No. 0134938 and a
%Sloan Research Fellowship.  The second author would also like to thank
%the Institut Mittag-Leffler (Djursholm, Sweden) for support during this work}
%FIX!!
The second author was partially supported by NSF Grant Number 0134938 and a
Sloan Research Fellowship.  The second author would also like to thank
the Institut Mittag-Leffler (Djursholm, Sweden) for support during this work.
\vspace{-3pt}

%%%%%%%%%%%%%%%%%%%%%%%%%%%%%%%%%%%%%%%%%%%%%%%%%%%%%%%%%%%%%%%%%%%%%

\section{Background on ring spectra}
In this section we give some of the basic constructions and properties
of ring spectra which will be used throughout the paper.
\vspace{-4pt}

\subsection{Postnikov sections}\label{2.1}
\vspace{-4pt}
Let $R$ be a connective, commutative ring spectrum.  For any $R$--module
$V$, we let $T_R(V)$ denote the tensor algebra on $V$.  For any
pointed simplicial set $K$, let $T_R(K)$ be shorthand for $T_R(R\Smash
\Sigma^\infty K)$.  Note that one has maps $T_R(\bdd{n})\ra T_R(\del{n})$,
and these are cofibrations of $R$--algebras.
\vspace{-4pt}

If $E$ is a cofibrant, connective
$R$--algebra and
\[ \xymatrix{ T_R(\bdd{n+1}) \ar[r]\ar[d] & E \ar[d] \\
T_R(\del{n+1}) \ar[r] & E'}
\]
is a pushout diagram of $R$--algebras, one verifies that $\pi_i(E)\ra
\pi_i(E')$ is an isomorphism for $i\leq n-1$ (see
\fullref{lem:smallBM}).  Here we need that $E$ is cofibrant as an
$R$--module (which follows from being cofibrant as an $R$--algebra, by
\cite[4.1(3)]{SS1})
%thm
to ensure that the pushout has the correct homotopy type---see the
proof of \fullref{lem:smallBM} for the details.
\vspace{-4pt}

If $E$ is a cofibrant, connective $R$--algebra, let $P_n(E)$ be the
result of applying the small object argument to $E$ with respect to
the set of maps $T_R(\bdd{i}) \ra T_R(\del{i})$ for all $i\geq n+2$
together with the generating trivial cofibrations for $R$--$\Alg$
from~\cite[5.4.3]{hss} or~\cite[4.1]{SS1}.  This is similar to the
functorial construction of a Postnikov section for differential graded
algebras given in~\cite[3.2]{tpwe}; see also \cite[5.1]{ms} for a
detailed description of functorial Postnikov sections for
symmetric spectra.  One 
checks that $P_n(E)$ is
fibrant, $\pi_i P_n(E)=0$ for $i>n$, and $\pi_i E\ra \pi_i P_n(E)$ is
an isomorphism for $i\leq n$.  Also, one has natural maps $P_{n+1}(E)
\ra P_n(E)$ which are compatible with $E\ra P_n(E)$ as $n$ varies.
\vspace{-4pt}

Finally, if $E$ is a connective $R$--algebra then we will write $P_n(E)$
as shorthand for $P_n(cE)$, where $cE\ra E$ is a fixed functorial
cofibrant-replacement for $E$.  Note that one does not have a map
$E\ra P_n(E)$ in this general case, only a zig-zag $E \bwe cE \ra
P_n(cE)=P_n(E)$.

\subsection{Pushouts of ring spectra}
For the result below we will need to use relative homotopy groups
$\pi_*(B,A)$ where $A\ra B$ is a map of spectra.
Note that we are not assuming that $A\ra B$ is a cofibration, as is
often done.  What we mean by $\pi_*(B,A)$ is therefore $\pi_*(W,A)$,
where we have functorially factored $A\ra B$ as $A\cof W \trfib B$.

At many places in the paper we will use the following Blakers--Massey
theorem for ring spectra: 

\begin{thm}
\label{thm:smallBM}
Let $R$ be a connective, commutative ring spectrum and let $m,n \geq 1$.
Suppose given a
homotopy pushout square of $R$--algebras
\[ \xymatrix{A \ar[r]\ar[d] & C \ar[d] \\ B \ar[r] & P
}
\]
in which the following conditions hold:
\begin{enumerate}[(i)]
\item $A$ is connective.
\item
$\pi_i(C,A)=0$ for $i < m $.
\item
$\pi_i(B,A)=0$ for $i< n $.
\end{enumerate}
Then $\pi_{i}(B,A)\ra \pi_{i}(P,C)$ is an isomorphism for $i
 < m + n -1 $ and a surjection for $i = m+n-1$. 
In particular, this means $\pi_i(P,C)= 0$ for $i < n$, which implies
$\pi_i C\ra \pi_i P$ is an isomorphism for $i< n -1 $ and a surjection
for $i=n - 1$.
\end{thm}
This is an important result; however, we have not been able to  
find a proof of it in the literature.  
See Goerss and Hopkins \cite[2.3.13]{GH} or Baues \cite[I.C.4.6]{baues}, though,
%thms
for full statements in related contexts.  For completeness we have included 
a proof of \fullref{thm:smallBM} in \fullref{se:pushout-proof}.

\subsection{Relative homotopy groups}
\label{se:relative}
Let $A\ra B$ be a map of $R$--algebras.  We claim that $\pi_*(B,A)$ is
in a natural way a bimodule over $\pi_*A$.  To explain this, it
suffices to assume that $A$ and $B$ are fibrant $R$--algebras, hence
fibrant as spectra.  Let $K\inc L$ be a cofibration of pointed
simplicial sets where $K$ is weakly equivalent to $\bdd{n}$ and $L$ is
contractible.  Then the relative homotopy group
$\pi_n(B,A)$ may be
described as equivalence classes of diagrams $\cD$ of the form
\[ \xymatrix{ \Sigma^\infty K  \ar[r] \ar[d] & A \ar[d] \\
\Sigma^\infty L \ar[r] & B,
}
\]
where two diagrams $\cD$ and $\cD'$ are equivalent if there is a diagram
\[ \xymatrix{ \Sigma^\infty [(K\times \del{1})/({*}\times \del{1})]
\ar[r] \ar[d] & A \ar[d] \\
\Sigma^\infty [(L\times \del{1})/({*}\times \del{1})] \ar[r] & B}
\]
which restricts to $\cD$ and $\cD'$ under the inclusions $\{0\}\inc
\del{1}$ and $\{1\}\inc \del{1}$, respectively.

\begin{remark}
To see why this description of $\pi_n(B,A)$ is valid, note that
the above equivalence class of diagrams is precisely $\pi_0$ of the
pullback of simplicial mapping spaces 
\[ \Map(\Sigma^\infty K,A) \times_{\Map(\Sigma^\infty K,B)}
\Map(\Sigma^\infty L, B).
\]
Our assumptions on $A$, $B$, and $K\ra L$ ensure that all the mapping
spaces are fibrant and that
$\Map(\Sigma^\infty L, B) \ra \Map(\Sigma^\infty K,B)$ is a fibration.
Since $L$ is contractible, $\Map(\Sigma^\infty L,B)$ is also
contractible.  So the above pullback is weakly equivalent to the
homotopy fiber of $\Map(\Sigma^\infty K,A) \ra \Map(\Sigma^\infty
K,B)$, which is what we want.
\end{remark}

Now suppose given an element $\alpha\in \pi_n(B,A)$, represented by a
diagram $\cD$ of the form
\[ \xymatrix{ \Sigma^\infty (\bdd{n})  \ar[r] \ar[d] & A \ar[d] \\
\Sigma^\infty (\del{n}) \ar[r] & B,
}
\]
as above.  Also assume given an element $\beta\in \pi_k A$ which is 
represented by a map $\Sigma^\infty(\bdd{k+1}) \ra A$.  Then one forms
the new diagram
\[ \xymatrixcolsep{1.2pc}\xymatrix{
 \Sigma^\infty (\bdd{n}\Smash \bdd{k+1})  \ar[r]^-\iso\ar[d] &
\Sigma^\infty(\bdd{n})\Smash \Sigma^\infty(\bdd{k+1})
\ar[r] \ar[d] & A \Smash A
\ar[d]  \ar@{=}[r] & A\Smash A \ar[r]^-\mu \ar[d] & A \ar[d] \\
\Sigma^\infty (\del{n}\Smash \bdd{k+1}) \ar[r]^-\iso &
\Sigma^\infty(\del{n})\Smash \Sigma^\infty(\bdd{k+1})
\ar[r] & B \Smash A \ar[r] &
B\Smash B \ar[r]^-\mu & B
}
\]
which represents a homotopy element of $\pi_{n+k}(B,A)$ (where we are
taking $K=\bdd{n}\Smash \bdd{k+1}$ and $L=\del{n}\Smash \bdd{k+1}$).
Here we
would have been slightly better off if we were working with
topological spaces rather than simplicial sets, as one can choose
homeomorphisms $S^{n-1}\Smash S^k  \iso S^{n+k-1}$ and $D^n\Smash S^k
\iso D^{n+k}$; but everything works out simplicially as
well, with only a little extra care.

We have just described a pairing $\pi_n(B,A)\times \pi_k A \ra
\pi_{n+k}(B,A)$, and one checks that this makes $\pi_n(B,A)$ into a
right module over $\pi_k A$.  A similar construction works for the
left module structure, and the verification that this gives a bimodule
is routine.

Note that in the long exact homotopy sequence of a pair, the connecting
homomorphism $\bd\co \pi_n(B,A) \ra \pi_{n-1}(A)$ is a map of
$\pi_*A$ bimodules.  This is because $\bd$ sends a homotopy element
represented by a diagram
\[ \xymatrix{ \Sigma^\infty (\bdd{n})  \ar[r] \ar[d] & A \ar[d] \\
\Sigma^\infty (\del{n}) \ar[r] & B
}
\]
to the element $\Sigma^\infty(\bdd{n}) \ra A$.

%%%%%%%%%%%%%%%%%%%%%%%%%%%%%%%%%%%%%%%%%%%%%%%%%%%%%%%

\section[k-Invariants for ring spectra: an outline]{$k$-Invariants for ring spectra: an outline}

In this section we give a basic outline of how $k$--invariants work for
ring spectra.  A $k$--invariant gives rise to a Postnikov extension,
and given a Postnikov extension we explain how to construct an associated
$k$--invariant.  These basic constructions will then be analyzed in a
more sophisticated way later in the paper.

\subsection{Extensions of ring spectra}
\label{se:exten-ring}
We continue to assume that $R$ is a connective, commutative ring
spectrum.
Fix an $n\geq 1$.  Let $C$ be a connective
$R$--algebra such that $P_{n-1}C\he C$ and
let $M$ be a $\pi_0(C)$--bimodule---that is, a
$(\pi_0 C)\tens_{\pi_0(R)} (\pi_0C)^{op}$--module.
As mentioned in the introduction, we wish to consider ring spectra 
$Y$ together with a map $Y\ra C$ such that
\begin{enumerate}[(i)]
\item $P_{n}Y\he Y$,
\item $P_{n-1}Y \ra P_{n-1}C$ is a weak equivalence,
\item $\pi_n(Y)\iso M$ as $\pi_0(Y)$--bimodules (where $M$ becomes a
$\pi_0(Y)$--bimodule via the isomorphism $\pi_0(Y)\ra \pi_0(C)$).
\end{enumerate}
The map $Y\ra C$ is called a \mdfn{Postnikov extension of $C$} of type
$(M,n)$.  Note that a particular choice of isomorphism $\pi_n(Y)\iso
M$ is not part of the data.

Let $\cM(C{+}(M,n))$ denote the category whose objects are such
Postnikov extensions; here a map from $X\ra C$ to $Y\ra C$ is a
weak equivalence $X\ra Y$ making the evident triangle commute.  We'll call
this category the \dfn{moduli space\/} of  Postnikov extensions of $C$ of type $(M,n)$.

If $\cC$ is a category, we'll write $\pi_0(\cC)$ for the connected
components of the nerve of $\cC$.  We wish to study
$\pi_0\cM(C{+}(M,n))$, as this will tell us how many ``homotopically different''
extensions of $C$ there are of type $(M,n)$.

\begin{defn}
A functor $F\co \cC \ra \cD$ between categories will be called a
\mdfn{weak equivalence} if it induces a weak equivalence on the
nerves.  The functor $F$ will be called a
\dfn{homotopy equivalence\/} if there is a functor $G\co \cD \ra \cC$
and zig-zags of natural transformations between $F\circ G$ and
$\id_\cD$, and between $G\circ F$ and $\id_\cC$.
\end{defn}

\begin{prop}
\label{pr:changebase1}
Suppose $C \ra C'$ is a weak equivalence of $R$--algebras.  There
is an evident functor $\rho\co \cM(C{+}(M,n)) \ra \cM(C'{+}(M,n))$
induced by composition, 
and this is a \weak 
equivalence.
\end{prop}

\begin{proof}
We will construct a homotopy inverse.
Given an object $X\ra C'$ of the category $\cM(C'{+}(M,n))$, functorially factor this map as
$X\trcof X_1 \fib C'$.  Letting $\phi(X)=C\times_{C'} X_1$,
then $\phi(X) \ra C$ is in $\cM(C{+}(M,n))$, using right properness
of $R$--$\Alg$.  So
this defines a functor $\phi\co \cM(C'{+}(M,n)) \ra \cM(C{+}(M,n))$.

It is simple to check that there is a zig-zag of natural weak
equivalences between the composite $\rho\circ \phi$ and the identity
map, and the same for the other composite $\phi\circ \rho$.
\end{proof}

By the above proposition, we can assume that $C$ is a cofibrant
$R$--algebra when studying $\cM(C{+}(M,n))$.  We will always be clear
when we are making this assumption, however.

\subsection{Bimodules}
By a $C$--bimodule we mean a left $(C\Smash_R C^{op})$--module.  As
remarked in the introduction, there is a $C$--bimodule $HM$ satisfying
$\pi_i(HM)=0$ for $i\neq 0$ and $\pi_0(HM)\iso M$ (as
$\pi_0(C)$--bimodules).  Moreover, a typical obstruction theory argument
shows that any two such bimodules are weakly equivalent.  We
abbreviate $HM$ to just $M$ in the rest of the paper.
\vspace{-3pt}

\begin{remark}
Note that the notion of $C$--bimodule depends on more than just the homotopy
type of $C$.  For if $C\ra C'$ is a weak equivalence of $R$--algebras,
the induced map $C\Smash_R C \ra C'\Smash_R C'$ need not be a weak
equivalence anymore.  For this reason we will sometimes have to assume
that $C$ is cofibrant as an $R$--module when dealing with bimodules.
\end{remark}
\vspace{-3pt}

\subsection{Extensions via pullbacks}
\label{se:exten-pullback}
\vspace{-3pt}
If $\cM$ is a model category and $X\in \cM$, let $\cM_{/X}$ denote the
usual overcategory whose objects are maps $Y\ra X$ in $\cM$.  Recall
that $\cM_{/X}$ inherits a model structure from $\cM$ in which a map from $Y\ra
X$ to $Y'\ra X$ is a cofibration (respectively fibration, weak equivalence)
if and only if the map $Y\ra Y'$ is a cofibration (respectively fibration,
weak equivalence) in $\cM$.
\vspace{-3pt}

We regard $C\Wedge\Sigma^{n+1}M$ as an object in $R$--$\smash{\Alg_{/C}}$ via the
projection $C\Wedge\Sigma^{n+1}M \ra C$.
Note that
$C\Wedge\Sigma^{n+1}M$ is actually a {\it pointed\/} object of
$R$--$\smash{\Alg_{/C}}$, since it comes equipped with the evident inclusion
$C\inc C\Wedge \smash{\Sigma^{n+1}}M$.
\vspace{-3pt}

Suppose given a homotopy class in $\ho(R$--$\smash{\Alg_{/C}})(C,C\Wedge
\Sigma^{n+1}M)$.  This can be represented by a map
\[ \alpha\co cC \ra f(C\Wedge\Sigma^{n+1}M) \]
where $cC$ is a cofibrant-replacement of $C$ and
$f(C\Wedge\Sigma^{n+1}M)$ is a fibrant-replace\-ment of
$C\Wedge\Sigma^{n+1}M$ in $R$--$\smash{\Alg_{/C}}$.  Consider the homotopy
fiber of $\alpha$ in $R$--$\smash{\Alg_{/C}}$.  This is the same as the homotopy
pullback of
\[ cC \ra f(C\Wedge \Sigma^{n+1}M) \la C \]
in $R$--$\Alg$.  
To be precise, 
to form this homotopy pullback we functorially factor the maps as
trivial cofibrations followed by fibrations
\[ cC \trcof (cC)' \fib f(C\Wedge\Sigma^{n+1}M),  \qquad \qquad
C \trcof C' \fib f(C\Wedge\Sigma^{n+1}M) \]
and then the homotopy pullback $Y$ is the pullback
\[\xymatrix{
Y  \ar@{.>}[r]\ar@{.>}[d] & C'\ar[d] \\
(cC)' \ar[r] & f(C\Wedge \Sigma^{n+1}M)
}
\]
in $R$--$\Alg$.  As pullbacks in $R$--$\Alg$ are the same as
pullbacks in ordinary spectra, it is easy to analyze the homotopy
groups of $Y$.  One sees immediately that $\pi_i((cC)',Y)=0$ for $i\neq
n+1$ and $\pi_{n+1}((cC)',Y)\iso M$.  So $\pi_i(Y)=0$ for $i>n$,
$\pi_i Y\ra \pi_i(cC)'$ is an isomorphism for $i<n$, and the map
$\bd\co\pi_{n+1}((cC)',Y) \ra \pi_n Y$ is an isomorphism.  By the
remarks at the end of \fullref{se:relative}, $\bd$ is an
isomorphism of $\pi_0(Y)$--bimodules.  Moreover, the map
$\pi_{n+1}((cC)',Y) \ra \pi_{n+1}(f(C\Wedge \Sigma^{n+1}M),C)$ is an
isomorphism of $\pi_0(Y)$--bimodules, and the codomain of this map is
clearly isomorphic to $M$ as a bimodule.  We have therefore shown that
$Y$ is a Postnikov extension of $(cC)'$ of type $(M,n)$.  As we have a
map $(cC)'\ra f(C\Wedge\Sigma^{n+1}M) \ra C$ which is a weak
equivalence, this is also a Postnikov extension of $C$ of type $(M,n)$.

The above remarks give us a function
\[ \PB\co \ho(R\text{--}\smash{\Alg_{/C}})(C,C\Wedge\Sigma^{n+1}M) \ra \pi_0\cM(C{+}(M,n)).
\]
It is clearly not injective, for the following reason.  An
automorphism $\sigma\co M\ra M$ of $\pi_0(C)$--bimodules induces an
automorphism of ring spectra $\sigma\co C\Wedge\Sigma^{n+1}M \ra
C\Wedge\Sigma^{n+1}M$.  If a given homotopy class $\alpha \in
\ho(R$--$\smash{\Alg_{/C}})(C,C\Wedge\Sigma^{n+1}M)$ is composed with this
$\sigma$, it gives rise to a weakly equivalent pullback.

Let $\Aut(M)$ be the group of $\pi_0(C)$--bimodule automorphisms of $M$.
One way to rephrase the above paragraph is to say that we have an
action of $\Aut(M)$
on the set of homotopy
classes we're considering, and we get an induced map
\[ \widetilde{\PB}\co  [\ho(R\text{--}\Alg_{/C})(C,C\Wedge\Sigma^{n+1}M)]/\Aut(M)
\ra \pi_0\cM(C{+}(M,n)).
\]
Our main goal in this paper is to show that this map is an
isomorphism.  Along the way, however, we will actually describe the
entire homotopy type of the space $\cM(C{+}(M,n))$ as opposed to just $\pi_0$.

\subsection[k-Invariants]{$k$-Invariants}
Let $f\co Y\ra C$ be a Postnikov extension of type $(M,n)$.  We
wish to show that it's in the image of $\widetilde{\PB}$.  We'll now
give a rough outline of how to go about this, which will then be
cleaned up in the later sections of the paper.

First of all, we can assume $Y$ is a cofibrant $R$--algebra (otherwise
we replace it with one).
Let $D$ be the homotopy pushout $C \amalg_Y^h C$ of 
$C \smash{\stackrel{f}{\vphantom{-}\smash{\lla}}} Y \smash{\stackrel{f}{\vphantom{-}\smash{\lra}}} C$ in
$R$--$\smash{\Alg_{/C}}$.  This means we factor $f$ as a cofibration followed
by  a trivial fibration
\[ Y \cof C' \trfib C \]
and we let $D$ be the pushout
\[ \xymatrix{
Y \ar@{   >->}[r]\ar@{   >->}[d] & C' \ar@{.>}[d] \\ 
C' \ar@{.>}[r] & D
}
\]
in $R$--$\Alg$.  Note that there is a map $D\ra C$ induced by the
universal property of pushouts.  Applying
\fullref{thm:smallBM}, we have that $\pi_i(C',Y) \ra
\pi_i(D,C')$ is an isomorphism for $i\leq 2n$.  It follows that
$\pi_i C' \ra \pi_iD$ is an isomorphism for $i\leq n$, and
$\pi_{n+1}(D,C')\iso M$ as $\pi_0(C')$--bimodules.  (For the bimodule
aspect of the last claim, one again uses the remarks in
\fullref{se:relative}).

Let $E=P_{n+1}D$ (and note that $D$ is cofibrant, so that we have a
natural map $D\ra P_{n+1}D$).  We will later show
that the map $C'\ra E$ is weakly equivalent in
$R$--$\smash{\Alg_{/C}}$ to the standard inclusion $C \inc C\Wedge\Sigma^{n+1}M$.
This can easily be done by an obstruction theory argument (see also
\fullref{re:EM}).  After
choosing such a weak equivalence, we have that the composite map
\[ C' \ra D \ra E \]
represents an element of $\Ho(R$--$\smash{\Alg_{/C}})(C,C\Wedge\Sigma^{n+1}M)$.
We will show that choosing a different weak equivalence only affects
this element up to the action of $\Aut(M)$, so that we have a
well-defined invariant in
\[ [\Ho(R\text{--}\Alg_{/C})(C,C\Wedge\Sigma^{n+1}M)]/\Aut(M).
\]
  Some unpleasant
checking is then required to verify that  we have indeed produced an
inverse to $\widetilde{\PB}$.  To organize this checking, it helps to
rephrase everything in terms of categories---this is what we do in the
next section.

%%%%%%%%%%%%%%%%%%%%%%%%%%%%%%%%%%%%%%%%%%%%%%%%%%%%%%%

\section[Categories of k-invariants and Eilenberg--Mac Lane objects]{\!\!Categories of $k\!$--invariants and Eilenberg--Mac Lane objects}

Let $R$, $C$, $M$, and $n$ be as in the previous section.
Our goal for the remainder of the paper is to analyze the homotopy
type of the moduli space $\cM(C{+}(M,n))$.  To do this we need to
introduce some auxiliary categories.

Define the category $\cE_C(M,n)$ of \mdfn{$C$--Eilenberg--Mac Lane
objects} of type $(M,n)$ as follows.
The objects of $\cE_C(M,n)$ are maps $B\ra E$ in $R$--$\smash{\Alg_{/C}}$ such that
\begin{enumerate}[(i)]
\item $B\ra C$ is a weak equivalence,
\item $B\ra E$ becomes a weak equivalence after applying $P_n$,
\item $P_{n+1}E\he E$,
\item $\pi_{n+1}(E)\iso M$ as $\pi_0(B)$--bimodules (where $M$ becomes
a $\pi_0(B)$--bimodule via the isomorphism $\pi_0(B)\ra \pi_0(C)$).
\end{enumerate}
A map from $[B\ra E]$ to $[B'\ra E']$ in this category consists of
weak equivalences $B\ra B'$ and $E\ra E'$ in $R$--$\smash{\Alg_{/C}}$ making the
evident square commute.

\begin{remark}
\label{re:EM}
Although it is not entirely obvious, we will see later that every
object of $\cE_C(M,n)$ is weakly equivalent to $C\inc
C\Wedge\Sigma^{n+1}M$.  This follows from \fullref{pr:prop2}(b)
below, which shows that $\cE_C(M,n)$ is connected.
\end{remark}

Likewise, we
define the category $\cK_C(M,n)$ of \mdfn{generalized
$k$--invariants for $C$ of type $(M,n)$}.
The objects of $\cK_C(M,n)$ are pairs of maps
$A\ra E \la B$ in the category $R$--$\smash{\Alg_{/C}}$ such that
\begin{enumerate}[(i)]
\item $A\ra C$ is a weak equivalence,
\item $B\ra E$ is an object in $\cE_C(M,n)$.
\end{enumerate}
A map from $A\ra E \la B$ to $A' \ra E' \la B'$ is a commutative diagram
\[\xymatrix{
A \ar[r] \ar[d]^\he & E \ar[d]^\he & B\ar[l]\ar[d]^\he \\
A'\ar[r] & E' & B'. \ar[l]
}
\]
in $R$--$\smash{\Alg_{/C}}$ in which all the vertical maps are weak equivalences.

Note that there is a forgetful functor $\cK_C(M,n)\ra
\cE_C(M,n)$ which forgets the object $A$ and the map $A\ra E$.

\begin{prop}
\label{pr:changebase2}
Suppose $C \ra C'$ is a weak equivalence of $R$--algebras.  Then there
are evident functors $\rho\co \cK_{C}(M,n) \ra \cK_{C'}(M,n)$
and $\rho\co \cE_C(M,n) \ra \cE_{C'}(M,n)$ induced by composition,
and both are \weak equivalences.
\end{prop}

\begin{proof}
We will prove the result for $\cK_C(M,N)$ and leave the $\cE_C(M,n)$
case to the reader.
Given $A\ra E \la B$ in
$\cK_{C'}(M,n)$, produce functorial factorizations
$A\trcof A' \fib C'$, $B\trcof B'\fib C'$, and $E\trcof E'\fib C'$.
So we have the diagram
\[\xymatrix{
A \ar[r] \ar[d]^\he & E \ar[d]^\he & B\ar[l]\ar[d]^\he \\
A'\ar[r] & E' & B'. \ar[l]
}
\]
Define $\phi\co \cK_{C'}(M,n) \ra \cK_{C}(M,n)$ by sending $A\ra E
\la B$ to the sequence of maps $A'\times_{C'} C \ra E'\times_{C'} C \la
B'\times_{C'} C$.  It follows by right-properness of $R$--$\Alg$ that this
sequence indeed lies in $\cK_{C}(M,n)$.

Just as in \fullref{pr:changebase1},
it is simple to produce a zig-zag of natural weak equivalences
between $\rho\circ \phi$ and the identity, and the same for
$\phi\circ \rho$.
\end{proof}

Our next goal is to show that the category $\cK_C(M,n)$
is weakly equivalent to the moduli space $\cM(C{+}(M,n))$.
First, observe that there is a functor $\PB_C\co \cK_C(M,n)
\ra \cM(C{+}(M,n))$ which sends $A\ra E \la B$ to its homotopy pullback
in $R$--$\smash{\Alg_{/C}}$.  
As in \fullref{se:exten-pullback},
whenever we talk about the ``homotopy
pullback'' $A \times_E^h B$ of a diagram $A\ra E \la B$ we mean the 
pullback of $A' \ra
E \la B'$ where we have functorially factored $A\ra E$ and $B\ra E$ as
trivial cofibrations followed by fibrations
\[ A \trcof A' \fib E \qquad \text{and}\qquad B \trcof B' \fib E.
\]
Note that there is a natural map from the
pullback of  $A\ra E \la B$ to its homotopy pullback.

To verify that the image of $\PB_C$ actually lands in $\cM(C{+}(M,n))$,
first recall that pullbacks of $R$--algebras are the same as pullbacks
of spectra.  This immediately verifies conditions (i) and (ii) in the
definition of $\cM(C{+}(M,n))$ (\fullref{se:exten-ring}).  For the
third condition, let $P=A'\times_E B'$.  Note that in the long exact
homotopy sequence of a pair the connecting homomorphism
$\pi_{n+1}(A',P)\ra \pi_n(P)$ is an isomorphism since $P_{n-1}A'\he
A'$, and in fact it is an isomorphism of $\pi_0(P)$--bimodules by the
discussion in \fullref{se:relative}.  But we also have an
isomorphism of $\pi_0(P)$--bimodules $\pi_{n+1}(A',P)\ra
\pi_{n+1}(E,B')$, as well as an isomorphism 
$\pi_{n+1}(E)\ra\pi_{n+1}(E,B')$ of $\pi_0(B')$--bimodules.  For the fact that these are
bimodule maps, we again refer to \fullref{se:relative}.  Since
$\pi_{n+1}(E,B')\iso \pi_{n+1}(E,B)\iso M$, we find that $\pi_n(P)$ is
isomorphic to $M$ as $\pi_0(P)$--bimodules.

\begin{prop}
\label{pr:prop1}
The functor $\PB_C\co \cK_C(M,n) \ra \cM(C{+}(M,n))$ is a weak equivalence.
\end{prop}

The proof is somewhat long, and will be given in \fullref{a.7}.
The basic idea is to try to construct a homotopy inverse functor
$k\co \cM(C{+}(M,n))\ra \cK_C(M,n)$.  This will be our
\dfn{generalized $k$--invariant\/}.  Given an $X\ra C$ in $\cM(C{+}(M,n))$
such that $X$ is cofibrant, construct a homotopy pushout
\[ \xymatrix{ X\ar[r]\ar[d] & P_{n-1}X \ar[d]^g \\
 P_{n-1}X \ar[r]^h & Z}
\]
and consider the maps
\[ P_{n-1}X\llra{h} P_{n+1}Z \llla{g} P_{n-1}X.
\]
\newpage

One can check that this gives an element $k(X)$ in
$\cK_{P_{n+1}(P_{n-1}C)}(M,n)$.  This
is \textit{almost\/} an element of $\cK_C(M,n)$ since $P_{n+1}(P_{n-1}C)\he C$. 
Some care is required in
getting around this small difference, and this is part of what is
accomplished in \fullref{a.7}.

The above proposition reduces the problem of studying $\cM(C{+}(M,n))$
to that of studying $\cK_C(M,n)$.  We do this by analyzing the
forgetful functor $\cK_C(M,n) \ra \cE_C(M,n)$.
We will prove the following in \fullref{a.11}:

\begin{prop}~
\label{pr:prop2}
\begin{enumerate}[(a)]
\item There is a homotopy fiber sequence of spaces
\[ \underline{R\text{--}\Alg_{\smash{/C}}}(C,C\Wedge\Sigma^{n+1}M) \ra |\cK_C(M,n)| \ra
|\cE_C(M,n)|, \]
where the first term denotes the homotopy function complex in the model
category $R$--$\smash{\Alg_{/C}}$.
\item There is a weak equivalence of spaces $|\cE_C(M,n)|\he
B\Aut(M)$ where $\Aut(M)$ is the group of automorphisms of $M$ as a
$\pi_0(C)$--bimodule.
\end{enumerate}
\end{prop}

Part (a) is basically a routine ``moduli space'' problem, of the type
considered in \cite[Section 2]{BDG}.  Part (b) involves similar
techniques but also
requires some careful manipulations of ring spectra.

%%%%%%%

\section[The moduli space of k-invariants]{The moduli space of $k$--invariants}\label{a.7}

In this section we will prove \fullref{pr:prop1}.  As
mentioned above, our first hope would be for a generalized
k-invariant functor $k\co \cM(C{+}(M,n)) \to \cK_C(M,n)$ to provide a
homotopy inverse for $\PB_C\co \cK_C(M,n) \to \cM(C{+}(M,n))$.  This
doesn't quite work out.  Instead, we restrict to the case where $C$ is
cofibrant and construct a functor $k\co \cM(C{+}(M,n)) \to
\cK_{C'}(M,n)$ where $C' = P_{n+1}(P_{n-1}C)$.  We then use this to
show that $\PB_C$ is a weak equivalence.

For the rest of this section we assume that $C$ is a cofibrant
$R$--algebra.  Throughout the following, let $C' = P_{n+1}(P_{n-1}C)$.
Note that since we are using the functor $P_{n+1}$ coming from
the small object argument described in \fullref{2.1}, the
map $P_{n-1} C \to P_{n+1}(P_{n-1}C)$ is only a weak equivalence
and not an isomorphism.   
Next we define $\cD(C,M,n)$, a category of diagrams which will be
useful in defining a generalized k-invariant functor $k\mc
\cM(C{+}(M,n)) \to \cK_{C'}(M,n)$.

The objects of $\cD(C,M,n)$ are the commutative diagrams $D$ of the
form
\[ \xymatrix{ H \ar[r]\ar[d] & B \ar[d] \\
 A \ar[r] & C}
\]
in which $A\ra C$ and $B\ra C$ are weak equivalences after applying
$P_{n-1}$, and $H\ra C$ lies in $\cM(C{+}(M,n))$.
The morphisms in $\cD(C,M,n)$ are the maps of commuting diagrams.
Using the weak equivalence $C \to P_{n-1}C \to P_{n+1}(P_{n-1}C) = C'$, we
will produce a diagram
\[ \xymatrix{
\cA(D) \ar[r]\ar[dr]& \cE(D) \ar[d] & \cB(D) \ar[l]\ar[dl] \\
& C'
}
\]
which is functorial in $\cD(C,M,n)$ and has the following properties:
\begin{enumerate}[(1)]
\item The diagram $\cA(D) \ra \cE(D) \la \cB(D)$ lies in
$\cK_{C'}(M,n)$, where we have used the weak equivalence $C\ra C'$
to make $M$ a bimodule over $\pi_0(C')$.
\item There is a natural zig-zag of weak equivalences in
$R$--$\Alg_{/C'}$ between  $H$ and the homotopy  pullback of $\cA(D) \ra
\cE(D) \la \cB(D)$.
\item If $A\fib E\bfib B$ is an object in $\cK_C(M,n)$ in which the
indicated maps are fibrations, and $D$ is the
diagram
\[\xymatrix{ A\times_E B \ar[r]\ar[d] & B\ar[d] \\
                 A\ar[r] & C
}
\]
then there is a natural zig-zag of weak equivalences between $A\ra E
\la B$ and $\cA(D) \ra \cE(D) \la \cB(D)$ in $R$--$\Alg_{/C'}$.
\item Suppose $D'$ is another diagram
\[ \xymatrix{H' \ar[r]\ar[d] & B'\ar[d] \\
A'\ar[r] & C}
\]
in $\cD(C,M,n)$, and assume there is a map of
diagrams $D\ra D'$ which is the identity on $C$, a weak equivalence on
$H\ra H'$ and on $B\ra B'$, and a weak equivalence after applying
$P_{n-1}$ to $A\ra A'$.  Then the induced maps $\cA(D)\ra \cA(D')$,
$\cE(D)\ra \cE(D')$, and $\cB(D)\ra \cB(D')$ are weak equivalences.
\end{enumerate}

In a moment we will explain how to construct $\cA(D)$, $\cE(D)$, and
$\cB(D)$, and we will verify the above properties.  But first
we show how this implies what we want.

\begin{proof}[Proof of \fullref{pr:prop1}]
Using \fullref{pr:changebase1} and \fullref{pr:changebase2}, it
suffices to analyze the case when $C$ is a cofibrant $R$--algebra.

Suppose given $X\ra C$ in $\cM(C{+}(M,n))$.  We define $k(X)\in \cK_{C'}(M,n)$
to be $\cA(D) \ra \cE(D) \la \cB(D)$ where $D$ is the diagram
\[\xymatrix{
X \ar[r]\ar[d]_{\id} & C \ar[d]^{\id} \\
X\ar[r] & C.
}
\]
We think of $k(X)$ as the \mdfn{generalized $k$--invariant of $X$}.
It gives a functor $k$ from $\cM(C{+}(M,n))$ to $\cK_{C'}(M,n)$.

We now have the following (noncommutative) diagram of functors:
\[ \xymatrix{
\cM(C{+}(M,n)) \ar[r]^{\rho_1}\ar[dr]^k & \cM(C'{+}(M,n)) \\
\cK_C(M,n)\ar[u]^{\PB_C}\ar[r]_{\rho_2} & \cK_{C'}(M,n) \ar[u]_{\PB_{C'}}
}
\]
The maps labelled $\rho_1$ and $\rho_2$ are known to be weak
equivalences, by \fullref{pr:changebase1} and
\fullref{pr:changebase2}.  If we can show that there is a zig-zag of
natural transformations between $\rho_1$ and the composite
$\PB_{C'}\circ k$, as well as between $\rho_2$ and the composite
$k\circ \PB_C$, it will follow that all maps in the diagram induce
isomorphisms on the homotopy groups of the nerves---so all the maps
will be weak equivalences.

Now, property (2) says precisely that there is a zig-zag of natural
weak equivalences between $\PB_{C'}\circ k$ and
$\rho_1$.  So we consider the other
composite.  Let $A\ra E \la B$ be an object in $\cK_C(M,n)$.  Functorially
factor $A\ra E$ and $B\ra E$ as
\[ A\trcof A' \fib E \qquad\text{and}\qquad B\trcof B'\fib E
\]
and let $H=A'\times_E B'$.  So $\PB_C(A\ra E\la B)=H$.
Let $D_1$, $D_2$, and $D_3$ be the following three squares:
\[\xymatrix{
H \ar[r]\ar[d] & C \ar[d] && H\ar[r]\ar[d] & C \ar[d] && H\ar[r]\ar[d] &
B'\ar[d] \\
H\ar[r] & C && A'\ar[r] & C && A' \ar[r] & C.
}
\]
Note that there are natural transformations $D_1\ra D_2$ and
$D_3\ra D_2$.  By property (4), we get the following chain of
equivalences:
\[ \xymatrix{
\cA(D_1) \ar[r] \ar[d]_{{\hbox{\footnotesize $\sim$}}} & \cE(D_1)\ar[d]_{\hbox{\footnotesize $\sim$}} &
\cB(D_1)\ar[l]\ar[d]_{\hbox{\footnotesize $\sim$}} \\
\cA(D_2)\ar[r] & \cE(D_2) & \cB(D_2) \ar[l] \\
\cA(D_3)\ar[u]^{\hbox{\footnotesize $\sim$}} \ar[r] & \cE(D_3)\ar[u]^{\hbox{\footnotesize $\sim$}} & \cB(D_3).\ar[l]\ar[u]^{\hbox{\footnotesize $\sim$}}
}
\]
Note that the top row is $k(H)=k(\PB_C(A\ra E \la B))$.
Using property (3), there is a natural zig-zag of weak equivalences
between the last row and the diagram $A' \ra E \la B'$, which in turn
is weakly equivalent to $A\ra E \la B$.  So we have established our
zig-zag of natural transformations between $k\circ \PB_C$ and $\rho_2$.
This completes the proof.
\end{proof}

Finally we are reduced to doing some actual work: we must
construct the functors $\cA$, $\cE$, and $\cB$.  Suppose given a
diagram $D$ in $\cD(C,M,n)$ of the form
\[ \xymatrix{ H \ar[r]\ar[d] & B \ar[d] \\
 A \ar[r] & C.}
\]
Recall that $A\ra C$ and $B\ra C$ are weak equivalences after applying
$P_{n-1}$, and $H\ra C$ lies in $\cM(C{+}(M,n))$.
Let $cX\we X$ be a functorial cofibrant-replacement in $R$--$\Alg$.
Consider the
composites $cH\ra P_{n-1}(cH) \ra P_{n-1}(cA)$ and $cH\ra P_{n-1}(cH) \ra
P_{n-1}(cB)$: functorially factor them as
\[ cH \cof SA \trfib P_{n-1}(cA) \qquad\text{and}\qquad cH\cof SB \trfib
P_{n-1}(cB).
\]
We obtain a diagram
\[ \xymatrix{
SA \ar[r]\ar[dr]_{\hbox{\footnotesize $\sim$}} & P_{n+1}[(SA)\amalg_{cH} (SB)] \ar[d]
& SB \ar[l]\ar[dl]^{\hbox{\footnotesize $\sim$}} \\
& P_{n+1}(P_{n-1}C).
}
\]
We let $\cA(D)=SA$, $\cE(D)=P_{n+1}[(SA)\amalg_{cH} (SB)]$, $\cB(D)=SB$, and
$C'=P_{n+1}(P_{n-1}C)$.

Notice the following:
\begin{itemize}
\item Property (4) follows immediately from our definitions.
\item There is a natural map from $cH$ into the pullback
$\cA(D)\times_{\cE(D)} \cB(D)$, and of course a natural map from the
pullback to the homotopy pullback.  This gives a natural zig-zag
\[ H \bwe cH \ra \cA(D)\times_{\cE(D)}^h \cB(D)
\]
in $R$--$\Alg_{/C'}$.
\end{itemize}

In order to check the remaining properties we will need the following
two lemmas.  

\begin{lemma}
\label{lem:bimodule}
Fix $n\geq 1$.
Let $W$ be a connective $R$--algebra satisfying $P_{n-1}W\he W$, and let
$M$ be a $\pi_0(W)$--bimodule.
Let $A \la X \ra B$ be maps in $R$--$\Alg_{/W}$ where $A\ra W$ and $B\ra
W$ are weak equivalences and $X\ra W$ is in $\cM(W{+}(M,n))$.  Let $P$
denote the homotopy pushout $A\amalg_X^h B$.  Then $\pi_{n+1}P$ is
isomorphic to $\pi_n X$ as a $\pi_0(X)$--bimodule.
\end{lemma}

\begin{proof}
Consider the map $f\co \pi_{n+1}(A,X) \ra \pi_{n+1}(P,B)$, which is
an isomorphism by \fullref{thm:smallBM}.  This is readily
seen to be a map of $\pi_0(X)$--bimodules, using the observations
from \fullref{se:relative}.  Here we regard $\pi_{n+1}(P,B)$
as a $\pi_0(X)$--bimodule via the map of $\pi_0(R)$--algebras
$\pi_0(X)\ra\pi_0(B)$.
The map $\pi_{n+1}(P) \ra \pi_{n+1}(P,B)$ is
a map of $\pi_0(B)$--bimodules (and hence $\pi_0(X)$--bimodules, by
restriction) which is an isomorphism by our assumptions on $B$.

Finally, the connecting homomorphism $\pi_{n+1}(A,X)\ra \pi_n(X)$ is a
map of $\pi_0(X)$--bimodules which is an isomorphism by our assumptions
on $X$ and $A$.
So we have established that
$\pi_{n+1}(P)\iso \pi_n(X)$ as $\pi_0(X)$--bimodules.
\end{proof}

\begin{lemma}\label{lem:main} Let $n \geq 1$ and let $W,M$ be as in the
previous lemma.
\begin{enumerate}[(a)]
\item Suppose that $A \fib E \bfib B$ is in $\cK_W(M,n)$ (where the
indicated maps are fibrations), and let
$H=A\times_E B$.  Then the induced map $A\amalg_H^h B \ra E$ becomes a weak
equivalence after applying $P_{n+1}$.
\item Suppose we are given cofibrations $A \bcof X \cofib B$ in $R$--$\Alg_{/W}$
such that both $A\ra W$ and $B\ra W$ are weak equivalences, and where
$P_n X\he X$.  Also assume both $P_{n-1}X\ra P_{n-1}A$ and
$P_{n-1}X\ra P_{n-1}B$ are weak equivalences.  Then the diagram $A\ra
P_{n+1}(A\amalg_X B) \la B$ lies in $\cK_{P_{n+1}W}(\pi_n X,n)$.
\item Again
suppose the given cofibrations $A \bcof X \cofib B$ in $R$--$\Alg_{/W}$ satisfy
the same hypotheses as in (b).  Let $X'$ be the homotopy pullback of
$A\ra P_{n+1}(A\amalg_X B) \la B$.  Then the induced map $X\ra X'$ is
a weak equivalence.
\end{enumerate}
\end{lemma}

\begin{proof}
The statements in (a) and (c) follow directly from
\fullref{thm:smallBM}.  One should note that pullbacks (and
homotopy pullbacks) in the category of $R$--algebras are the same as
those in the category of symmetric spectra.

For the statement in (b), one uses \fullref{thm:smallBM}
together with \fullref{lem:bimodule} above.
\end{proof}

\begin{proof}[Proof of Properties (1)--(3)]
Properties (1) and (2) follow directly from \fullref{lem:main} parts
(b) and (c), respectively.  So we turn to property (3).

Let $A\fib E \bfib B$ be an object in $\cK_C(M,n)$, and let $H$ be the
pullback $A \times_E B$.  Let $\cD$ be the diagram
\[ \xymatrix{ H \ar[r]\ar[d] & B \ar[d] \\
 A \ar[r] & C.}
\]
Functorially factor the maps $cH\ra cA$ and $cH\ra cB$ as $cH\cof S'A
\trfib cA$ and $cH\cof S'B \trfib cB$.  Note that one gets induced maps
$S'A \ra SA$ and $S'B\ra SB$ (where $SA$ and $SB$ appeared in our
construction of $\cA(D)$, etc.), and these maps are weak equivalences.
Let $\cE'(D)=P_{n+1}(S'A \amalg_{cH} S'B)$, so that there is an induced
map $\cE'(D)\ra \cE(D)$.

Notice that we have a map $S'A\amalg_{cH} S'B \ra cA\amalg_{cH} cB \ra
cE$, and therefore get an induced map $\cE'(D) \ra P_{n+1}(cE)$.  This
is a weak equivalence by \fullref{lem:main}(a).
Now we have the following:
\[\xymatrix{
 SA \ar[r] & \cE(D) & SB \ar[l]\\
 S'A \ar[r]\ar[u]^{\hbox{\footnotesize $\sim$}}\ar[d]_{\hbox{\footnotesize $\sim$}}
     & \cE'(D)\ar[u]^{\hbox{\footnotesize $\sim$}}\ar[d]_{\hbox{\footnotesize $\sim$}} & S'B\ar[u]^{\hbox{\footnotesize $\sim$}}\ar[d]_{\hbox{\footnotesize $\sim$}} \ar[l] \\
cA \ar[r] & P_{n+1}(cE) & cB \ar[l] \\
 cA\ar[r]\ar[u]^{\hbox{\footnotesize $\sim$}} \ar[d]_{\hbox{\footnotesize $\sim$}}
& cE\ar[u]^{\hbox{\footnotesize $\sim$}}\ar[d]_{\hbox{\footnotesize $\sim$}} & cB\ar[l]\ar[u]^{\hbox{\footnotesize $\sim$}}\ar[d]_{\hbox{\footnotesize $\sim$}} \\
  A\ar[r] & E & B.\ar[l] \\
}
\]
So we have obtained a natural zig-zag of weak equivalences
between the diagrams $\cA(D) \ra \cE(D) \la \cB(D)$ and $A\ra E \la B$.
\end{proof}

%%%%%%%%%%%%%%%%%%%%%%%%%%%%%%%%%%%%%%%%%%%%%%%%%%%%%%%%%%%%%%%%%%%

\section{Nonunital algebras}\label{sec:nonu}

Before proceeding further with our main results we need to develop a
little machinery.  This concerns nonunital $C$--algebras and their
relations to $C$--bimodules.  We first discuss a model structure
on nonunital $C$--algebras.  Then we define an ``indecomposables''
functor from nonunital $C$--algebras to $C$--bimodules and study
its interaction with Postnikov stages.  All of this basic machinery
has been heavily influenced by Basterra \cite{B}.

Let $R$ be a commutative $S$--algebra and $C$ an $R$--algebra.
Define a \mdfn{nonunital $C$--algebra} to be a
nonunital monoid in the category of $C$--bimodules, that is,
an algebra over the monad
\[ \tT(M) = M \amalg (M \sm_C M) \amalg (M\sm_C M \sm_C M) \amalg
\cdots
\]
in $\Cbimod$.  Let $\NonU_C$ denote the category of nonunital
$C$--algebras.

A map of nonunital $C$--algebras is defined to be a {\em fibration\/} or
a {\em weak equivalence\/} if the underlying map in $C\bimod$ (or
$R\lMod$) is a fibration or a weak equivalence.  A map is then a {\em
cofibration\/} if it has the left lifting property with respect to all
trivial fibrations.  Below we will use~\cite{SS1} to verify that this
gives a model structure on $\NonU_C$.

First recall that $\Cbimod$ is just another name for the category
$C \sm_R C^{op}\lMod$.  The model structure on $R$--modules lifts to a
model structure on $C$--bimodules \cite[4.1]{SS1}.  Let $\CC\co R\lmod
\to \Cbimod$ be the left adjoint to the forgetful functor $\Cbimod \to
R\lmod$.  If $j\mc K \to L$ is a generating trivial cofibration in $R\lMod$, then the generating trivial cofibrations in $\Cbimod$ are of the
form $ \CC(j)\mc C \sm_R K \sm_R C \to C \sm_R L \sm_R C $.

\begin{thm}\label{thm-nonu}
The above notions of cofibration, fibration, and weak equivalence form
a cofibrantly generated model category structure on $\NonU_C$.  The
generating cofibrations and trivial cofibrations are of the form
$\tT(\CC(K)) \to \tT(\CC(L))$ where $K\to L$ is a generating
cofibration or trivial cofibration in $R\lmod$.
\end{thm}

\begin{proof}
To establish the model structure on $\NonU_C$ we modify the arguments
for unital monoids in~\cite[6.2]{SS1}.  The argument in \cite{SS1} is
mostly formal except for one key step.  For us, this step is to show that
given
a generating trivial cofibration $K \to L$ in $\Cbimod$,
the pushout in $\NonU_C$ of the diagram
\[ \xymatrix{ \tT(K)  \ar[r] \ar[d] & \tT(L)  \\
X &   } \]
is the colimit $P = \colim P_n$ in $C\bimod$ of a sequence
\[ X = P_0 \to P_1 \to \cdots \to P_n \to \cdots \]
where $P_n$ is obtained from $P_{n-1}$ by a pushout in $\Cbimod$.  We
then show that the monoid axiom implies that $X=P_0 \to P$ is a weak
equivalence.  From this it follows directly from~\cite[2.3(1)]{SS1} that the
given model structure exists on $\NonU_C$.

To construct $P_n$ from $P_{n-1}$,
we replace the functor $W$ in~\cite[6.2]{SS1} by a functor
\[ W\co {\mathcal P}(\underline n)  \to C\bimod \]
where ${\mathcal P}(\underline n)$ is the poset category of subsets of
$\underline n = \{1, 2, \cdots, n\}$ and inclusions.  For $S \subseteq
\underline n$, define
\[ W(S) = (C\Wedge X) \sm_C B_1 \sm_C (C\Wedge X) \sm_C B_2 \sm_C
\cdots
\sm_C B_n \sm_C (C\Wedge X)\]
\[ B_i \ =
\begin{cases} K & \text{if $i\not\in S$} \\
L & \text{if $i\in S$.}
\end{cases} \leqno{\hbox{with}}
\]
Let $Q_{n}$ be the colimit of $W(S)$ over ${\mathcal P}(\underline n)
- \underline n$ (that is, the proper subsets of $\underline n$).  As
in~\cite[6.2]{SS1}, one has maps $Q_n \to P_{n-1}$ and $Q_n \to
W(\underline n)$ and defines $P_n$ as the following
pushout in $C\bimod$
\[ \xymatrix{ Q_n  \ar[r] \ar[d] & W(\underline n)\ar[d]  \\
P_{n-1} \ar[r] & P_n  } \]
Set $P = \colim P_n$, the colimit in $C\bimod$.  Arguments analogous to
those in~\cite[6.2]{SS1} show that $P$
is naturally a nonunital $C$--algebra and has the universal property of
the pushout of $X \la \tT(K) \ra \tT(L)$ in nonunital $C$--algebras.

We next show that the monoid axiom for $R\lmod$ implies that each map
$P_{n-1} \to P_n$ is a weak equivalence whenever $K \to L$ is a generating
trivial cofibration in $C\bimod$.  These generating trivial cofibrations
are of the form $\CC(K') \to \CC(L')$ where $K' \to L'$ is a generating
trivial cofibration in $R\lmod$.

Since pushouts in $C\bimod$ are created in $R\lmod$ which is
symmetric monoidal, we consider the pushouts defining $P_n$ in
$R\lmod$.  Replacing $K\to L$ by $\CC(K')\to \CC(L')$, we see that
$Q_{n} \to W(\underline n)$ is isomorphic to
\[
{Q'_n} \sm_R
(C\Wedge X)^{\sm_R (n)} \to (L')^{\sm_R n} \sm_R (C\Wedge X)^{\sm_R
(n)}
\]
where
${Q'_n} \to (L')^{\sm_R n}$ is the $n$--fold box product of $K' \to
L'$.  The pushout product axiom implies that
${Q'_n} \to (L')^{\sm_R n}$ is a trivial cofibration.  The monoid
axiom then implies that the pushouts $P_{n-1} \to P_n$ are weak
equivalences and $X= P_0 \to \colim P_n$ is a weak equivalence.
\end{proof}

There is a functor
$K\co \NonU_C \to (C\ovcat R$--$\Alg\ovcat C)$
which takes a nonunital algebra $N$ to $K(N)= (C
\smash{\stackrel{\eta~}{\vphantom{-}\smash{\rightarrow}}} C \vee N \smash{\stackrel{~\pi}{\vphantom{-}\smash{\lla}}} C)$ where $\eta$ and
$\pi$ are the obvious inclusion and projection.  This has a right
adjoint
\[
I\co (C\ovcat R\text{--}\Alg\ovcat C) \to \NonU_C
\]
called the
{\em augmentation ideal\/} functor.  The functor $I$
sends $(C\ra X \ra C)$ to the fiber of $X\ra C$.

\begin{prop}\label{prop.KI}
The functors $(K,I)$ form a Quillen equivalence
\[ K\co\NonU_C \we (C\ovcat R\text{--}\Alg\ovcat C).
\]
\end{prop}

\begin{proof}
The same statement for commutative ring spectra is
proved in Basterra \cite[2.2]{B}.
% Prop.
The proof works verbatim in the noncommutative
case; see also Basterra and Mandell \cite[Theorem 8.6]{BM} for a vast generalization.
\end{proof}

\subsection{Indecomposables}
Our next task is to compare nonunital $C$--algebras with $C$--bimodules.
The {\em indecomposables\/} functor $Q\mc \NonU_C \to \Cbimod$
takes $X$ in $\NonU_C$ to the pushout of $* \from X\sm X \to X$.
Its right adjoint $Z\mc \Cbimod \to \NonU_C$ sends a bimodule
to itself equipped with the zero product.

\begin{prop}\label{prop.QZ}
The functors $Q$ and $Z$ form a Quillen pair:
\[ \xymatrix{
\NonU_C \ar@<0.5ex>[r]^Q & C\bimod. \ar@<0.5ex>[l]^Z
}
\]
\end{prop}

\begin{proof}
The functor $Z$ obviously preserves fibrations and trivial
fibrations.  Again, see~\cite[3.1]{B} and
\cite[Proposition 8.7]{BM} for similar statements.
\end{proof}

If $N$ is a $C$--bimodule, it is easy to see that $Q(ZN)\iso N$.  When
we consider the derived functors $\barQ(\barZ N)$ the situation
changes, however.  Here we must take a cofibrant replacement of $ZN$
before applying $Q$.  It turns out that $\barQ(\barZ N)$ typically has
nonzero homotopy groups in infinitely many dimensions, even if $N$ did
not.

If $n\geq 1$ and $N$ has no homotopy groups in dimensions smaller than
$n$, then the same turns out to be true for $\barQ(\barZ N)$.
Moreover, the $n$--th homotopy group of $\barQ(\barZ N)$ is easy to
analyze, and it is the same as that of $N$.  This is the content of
\fullref{pr:Qanalysis} below.  We remark that this proposition
bears some relation to \cite[8.2]{B}. %lemma

Before stating the proposition we need a couple of pieces of new
notation.  We'll use $c$ and $f$ to denote cofibrant- and
fibrant-replacement functors in a model category, and we leave it to
the reader to decide from context which model category the
replacements are taking place in.  In the statement of the proposition
below, for instance, the $c$ is being applied in $\NonU_C$ and the
$f$'s are being applied in $C\bimod$.

Also, note that the $P_n$'s in the statement of the
proposition refer to Postnikov sections in the category of
$C$--bimodules.  These are constructed analogously to those for ring
spectra, but here one forms pushouts with respect to the maps
$(C\Smash_R C^{op})\Smash \bdd{i} \ra (C\Smash_R C^{op})\Smash
\del{i}$ for $i > n+1$.

\begin{prop}
\label{pr:Qanalysis}
Fix $n\geq 1$, and
let $N$ be a $C$--bimodule such that $\pi_i(N)=0$ for $i<n$.  There is
a natural weak equivalence of $C$--bimodules
\[ P_n[QcZ(fN)] \ra P_n(fN).
\]
\end{prop}

\begin{proof}
First note that there are natural maps $QcZ(fN)\ra QZ(fN) \ra fN$.
Applying $P_n$ to this composite gives the map in the statement of the
proposition.  Call this map $g$.

Recall that $\tT\co C\bimod \ra \NonU_C$ is the left adjoint of the
forgetful functor.
Note that this is a left Quillen functor, and that
there are natural isomorphisms $Q\smash{\tT}(W)\iso W$ by an
easy adjointness argument.

If $K$ is a spectrum, write $FK$ as shorthand for $(C\Smash_R
C^{op})\Smash K$.  This is the free $C$--bimodule generated by $K$.

Write $\PP_n$ for the Postnikov section functor in the category
$\NonU_C$.
The construction is the same as for ring spectra, but in this case we
form pushouts with respect to the maps $\smash{\tT}(F\bdd{i})\ra
\smash{\tT}(F\del{i})$ for $i>n+1$.   
Using that $Q$ is a left adjoint and therefore preserves pushouts and
colimits, 
and that $Q\smash{\tT}(W) \iso W$, 
one can show that $Q(\PP_nX)$ is obtained from $QX$ by forming 
pushouts with respect to $F(\bdd{i}) \to F(\del{i})$ for $i > n+1$.    
Given the description above of Postnikov sections
for $C$--bimodules, this implies that the map $QX\ra Q(\PP_n X)$ induces a
weak equivalence
\begin{myequation}
\label{eq:eq1}
 P_n[QX]\ra P_n[Q\PP_n X].
\end{myequation}
Without loss of generality we may assume that $N$ is a cofibrant
$C$--bimodule.
Consider the natural map of nonunital algebras $\smash{\tT}(N) \ra
Z(N)$, adjoint to the isomorphism $N\ra UZ(N)$ where $U\co \NonU_C
\ra C\bimod$ is the forgetful functor.
Using our hypothesis on $N$, one finds that $N\Smash_C N \sm_C \cdots \sm_C N$
doesn't have any homotopy groups in dimensions less than $2n$ as long as 
there are at least two smash factors of $N$.  It follows immediately that
the induced map $\PP_n[\tT(N)] \ra \PP_n(ZN)$ is a weak equivalence.
\vspace{1pt}

Consider the trivial fibration $cZ(fN) \ra Z(fN)$.  Since $\tT(N)$ is
cofibrant (since $N$ is), the map $\tT(N)\ra Z(N) \ra Z(fN)$ lifts to
a map $\tT(N) \ra cZ(fN)$.  This becomes a weak equivalence after
applying $\PP_n$, by the previous paragraph.

In the square
\[ \xymatrix{
P_nQ[\tT N] \ar[r]\ar[d] & P_n Q [\PP_n (\tT N)] \ar[d] \\
P_n Q [cZ(fN)] \ar[r] & P_n Q [\PP_n (cZ(fN))]
}
\]
the two horizontal maps are weak equivalences by \eqref{eq:eq1}.  The
previous paragraph shows that the right vertical map is a weak
equivalence, so the left vertical map is as well.  Thus, we have an
equivalence
\[ P_n(N)\iso P_n[Q\tT(N)] \ra P_nQ[cZ(fN)].
\]
It is easy to use the adjoint functors to see that the composite of
this map with our map $g\co P_nQ[cZ(fN)] \ra P_n(fN)$ is the map
$P_nN \ra P_n(fN)$ induced by $N\ra fN$.  Since this composite is a
weak equivalence, so is the map $g$.
\end{proof}

\begin{prop}
\label{pr:Qanalysis2}
Fix $n\geq 1$, and
let $N$ be a $C$--bimodule such that $\pi_i(N)=0$ for $i\neq n$.
If $X \in \NonU_C$ is weakly equivalent to $Z(N)$,
then the natural map $cX \ra ZP_{n}[QcX]$ is a weak equivalence.
\end{prop}

\begin{proof}
The natural map in the statement is the composite $\eta\co cX \ra
ZQ[cX] \ra ZP_n[QcX]$.  Since there will necessarily be a weak
equivalence $cX \ra Z(fN)$, it suffices to check that $\eta$ is a weak
equivalence when $X=Z(fN)$.  In this case we consider the diagram
\[ \xymatrix{
cZ(fN) \ar[dr]_{\hbox{\footnotesize $\sim$}}\ar[r] & ZQ[cZ(fN)] \ar[r]\ar[d] & ZP_n Q[cZ(fN)]
\ar[d]  \\
& Z(fN) \ar[r] & ZP_n(fN).
}
\]
The composite of the top horizontal maps is $\eta$, the
vertical maps come from the counit of the $(Q,Z)$ adjunction,  and 
the diagram is readily checked to commute.  
The bottom horizontal map is a weak equivalence because
$fN\ra P_n(fN)$ is a weak equivalence, and $Z$ preserves all weak
equivalences.  Finally, we know by the preceding proposition that the
right vertical map is a weak equivalence, so $\eta$ is a weak
equivalence as well.
\end{proof}
\vspace{-5pt}

%%%%%%%%%%%%%%%%%%%%%%%%%%%%%%%%%%%%%%%%%%%%%%%%%%%%%%%%%%%%%%%%%%%%%%%

\section{The moduli space of Eilenberg--Mac Lane objects}\label{a.11}
\vspace{-5pt}
Recall where we are at this point in the paper.  We have completed the
proof of \fullref{pr:prop1}, and our next goal is to
prove \fullref{pr:prop2}.
\vspace{-5pt}

\subsection{General moduli space technology}
\vspace{-5pt}
If $\cC$ is a model category and $X$ is an object of $\cC$,  let
$\cM_\cC(X)$ denote the category consisting of all objects weakly
equivalent to $X$, where the maps are weak equivalences.  This is
called the \mdfn{Dwyer--Kan classification space} of $X$, or the
\mdfn{moduli space} of $X$.  It is a theorem of Dwyer--Kan
\cite[2.3]{classifi}
%prop
that
$|\cM_\cC(X)| \he B\hAut(X)$ where $\hAut(X)$ denotes the simplicial
monoid of homotopy automorphisms of $X$.
This is simply the
subcomplex of the homotopy function complex $\mC(X,X)$ consisting of
all path components which are invertible in the monoid $\pi_0\mC(X,X)$.
\vspace{-5pt}

If $X$ and $Y$ are two objects of $\cC$ then we'll write
$\ccH_{\cC}(X,Y)$ for the category consisting of diagrams
\[ \xymatrix{
X  & U \ar[l]_{{\hbox{\footnotesize $\sim$}}}\ar[r] & V & Y\ar[l]_{{\hbox{\footnotesize $\sim$}}}
}
\]
where the indicated maps are weak equivalences.  A morphism in this
category is a natural weak equivalence between diagrams which is the
identity on $X$ and $Y$.
It is another result of Dwyer--Kan that one has a natural zig-zag of weak
equivalences between $\ccH_{\cC}(X,Y)$ and the homotopy function
complex $\underline{\cC}(X,Y)$.  This follows from \cite[6.2(i),
8.4]{calculating}.
%both props
\vspace{-5pt}

When $Y$ is fibrant one may consider a simpler moduli space:
let $\ccH_{\cC}(X,Y)^f$ be the category whose objects are diagrams
\[
\xymatrix{
X  & U \ar[l]_{{\hbox{\footnotesize $\sim$}}}\ar[r] & Y.
}
\]
A map in this category is again a weak equivalence of diagrams which
is the identity on $X$ and $Y$.
There is an inclusion functor $\ccH_{\cC}(X,Y)^f \ra
\ccH_{\cC}(X,Y)$, and it is stated in \cite[2.7]{BDG}
%remark
that this is a weak
equivalence when $Y$ is fibrant.  For a proof, see \cite{dugnew}.
\vspace{-5pt}

\begin{remark}
All of the Dwyer--Kan theorems we mentioned above were actually proven
only for {\it simplicial\/} model categories.  It is not obvious
whether the category $R$--$\Alg$ is simplicial, though.  Using the main
result of \cite{D1}, however, the Dwyer--Kan results can be immediately
extended to all combinatorial model categories.  All of the model
categories considered in the present paper are combinatorial, so we
will freely make use of this technology.
\end{remark}

\subsection{Applications to ring spectra}
\label{se:theta}
Note that there is a functor
\[
\Theta \mc\ccH_{R\text{--}\Alg_{/C}}(C,C\Wedge\Sigma^{n+1}M) \ra \cK_C(M,n)
\]
which sends a
diagram
\[
\xymatrix{
C & U \ar[l]_-{{\hbox{\footnotesize $\sim$}}}\ar[r] & V & C\Wedge\Sigma^{n+1}M\ar[l]_-{{\hbox{\footnotesize $\sim$}}}
}
\]
in $R$--$\smash{\Alg_{/C}}$ to the object of $\cK_C(M,n)$ represented by
\[ U \ra V \la C \]
(where the second map is the composite $C\inc C\Wedge\Sigma^{n+1}M \we
V$).

The following result is very similar to \cite[2.11]{BDG}.
%Remark

\begin{lemma}
\label{lem:lem2a}
The sequence of maps
\[ \ccH_{R\text{--}\Alg_{/C}}(C,C\Wedge\Sigma^{n+1}M)
\ra \cK_C(M,n) \ra
\cM_{R\text{--}\Alg_{/C}}(C) \times \cE_C(M,n)
\]
is a homotopy fiber sequence of simplicial sets.
\end{lemma}

Note that $\cM_{R\text{--}\Alg_{/C}}(C)$ is contractible, as $C$ is a terminal
object of this category.  So this lemma, together with the
identification of $\ccH_{R\text{--}\Alg_{/C}}(X,Y)$ with the homotopy function
space $\underline{R\text{--}\Alg_{\smash{/C}}}(X,Y)$, yields \fullref{pr:prop2}(a).  

\begin{proof}
The second map is 
 $F\co \cK_C(M,n) \ra \cM_{R\text{--}\Alg_{/C}}(C) \times \cE_C(M,n)$ which
sends an object $A\ra E \la B$ to the pair consisting of
$A$ and $E\la B$.
The proof of this lemma will be
an application of Quillen's Theorem B---however, a little care is
required.

Let $\cK_C^f(M,n)$ denote the full subcategory of $\cK_C(M,n)$
consisting of objects $A\ra E \la B$ where the maps $E\ra C$ and $B\ra
C$ are fibrations.  Let $\smash{\cE_C^f(M,n)}$ denote the analogous
subcategory of $\cE_C(M,n)$.  The inclusions $\smash{\cK_C^f(M,n)}\inc
\cK_C(M,n)$ and $\smash{\cE_C^f(M,n)\inc \cE_C(M,n)}$ are readily checked to
be homotopy equivalences.  Let $$\widetilde{F}\co \smash{\cK_C^f(M,n) \ra
\cE_C^f(M,n)}$$ be the restriction of the functor $F$.

To apply Quillen's Theorem B \cite{quillen}, we are required to check that for every
map $[E'\la B'] \ra [E''\la B'']$ in $\smash{\cE_C^f(M,n)}$ the induced map of
overcategories
\[ (\widetilde{F}\ovcat [E'\la B']) \ra (\widetilde{F}\ovcat [E''\ra B''])
\]
is a weak equivalence.
There is a functor
\[ \phi'\co \ccH_{R\text{--}\Alg_{/C}}(C,E')^f \ra (\widetilde{F}\ovcat
[E'\la B'])
\]
sending a diagram $C \bwe A \ra E'$ to the pair consisting
of the object $A\ra E' \la B'$ in $\smash{\cK_C^f(M,n)}$ together with the
identity map from $\widetilde{F}(A\ra E' \la B')$ to $[E'\la B']$.
This functor $\phi'$ is readily checked to be a homotopy equivalence.
Similarly, one has
\[ \phi''\co \ccH_{R\text{--}\Alg_{/C}}(C,E'')^f \ra (\widetilde{F}\ovcat
[E''\la B''])
\]
which is a homotopy equivalence by the same argument.  The map
\[ \ccH_{R\text{--}\Alg_{/C}}(C,E')^f \ra
\ccH_{R\text{--}\Alg_{/C}}(C,E'')^f
\]
is a weak equivalence because it is naturally equivalent to the map of
function complexes $\underline{R\text{--}\Alg_{\smash{/C}}}(C,E') \ra
\underline{R\text{--}\Alg_{\smash{/C}}}(C,E'')$ (which is itself a weak equivalence
because $E' \ra E''$ is a weak equivalence).  So we have established
that our map of overcategories is a weak equivalence.

Quillen's Theorem B now tells us that the sequence
\[
(\widetilde{F}\ovcat [f(C\Wedge\Sigma^{n+1}M) \la C]) \ra \cK_C^f(M,n)
\llra{\widetilde{F}} \cE_C^f(M,n)
\]
is a homotopy fiber sequence, where the basepoint in the base space is
taken to be the object $[f(C\Wedge\Sigma^{n+1}M) \la C]$.
We have already remarked that the
overcategory appearing here is homotopy equivalent to
the moduli category
$$\ccH_{R\text{--}\Alg_{/C}}(C,f(C\Wedge\Sigma^{n+1}M))^f,$$ and that the
inclusions $\cK_C^f(M,n)\inc
\cK_C(M,n)$ and $\cE_C^f(M,n)\inc \cE_C(M,n)$ are homotopy equivalences.
So to complete the proof of the lemma it suffices to note two things.
First, we have the commutative square
\[ \xymatrix{
\ccH_{R\text{--}\Alg_{/C}}(C,f(C\Wedge\Sigma^{n+1}M)) \ar[r] & \cK_C(M,n)
\\
\ccH_{R\text{--}\Alg_{/C}}(C,f(C\Wedge\Sigma^{n+1}M))^f \ar[r]\ar[u]_{\hbox{\footnotesize $\sim$}} &
\cK^f_C(M,n) \ar[u]_{{\hbox{\footnotesize $\sim$}}}
}
\]
where both vertical maps are the inclusions and the horizontal
maps are induced by the functor $\Theta$
defined in \fullref{se:theta}. 
Second, one has
a commutative triangle
\[ \xymatrix{
\ccH_{R\text{--}\Alg_{/C}}(C,f(C\Wedge\Sigma^{n+1}M)) \ar[r]\ar[d]_{{\hbox{\footnotesize $\sim$}}} &
\cK_C(M,n) \\
\ccH_{R\text{--}\Alg_{/C}}(C,C\Wedge\Sigma^{n+1}M) \ar[ur]
}
\]
where, again, all the maps are the obvious ones.
\end{proof}

Our next goal is to
identify $|\cE_C(M,n)|$ as in \fullref{pr:prop2}(b).  We first
show that for $C$ a cofibrant $R$--algebra, $\cE_C(M,n)$ is 
equivalent to a moduli space in the category
$\NonU_C$ of nonunital $C$--algebras.  This then reduces further to a
moduli space in the category of $C$--bimodules.  Computations in the
category of $C$--bimodules are relatively simple, so that it is not
hard to identify the homotopy type of the moduli space in
$C$--bimodules as $B\Aut(M)$.  Finally, we remove the cofibrancy
condition in \fullref{7.8}.

\begin{lemma}
\label{lem:lem2b}
Assume $C$ is a cofibrant $R$--algebra.
There are weak equivalences of categories
\[ \cE_C(M,n) \he \cM_{(C\ovcat R\text{--}\Alg \ovcat C)}(C\Wedge
\Sigma^{n+1}M) \he \cM_{\NonU_C}(\Sigma^{n+1}M).\]
\end{lemma}

\begin{proof}
We write $\cE$ for $\cE_C(M,n)$.
The argument proceeds in several steps.  First, let $\cE'$ be the
full subcategory of $\cE$ whose objects $B\ra E$ have $B=C$ (and the
map $B\ra C$ the identity).  Let $\cE''$ be the
full subcategory of $\cE$ whose objects are maps $B\cof E$ in $R$--$\smash{\Alg_{/C}}$
in which both $B$ and $E$ are cofibrant $R$--algebras and $B\ra E$ is a
cofibration.  Finally, let $\cE'''$ be the full subcategory of $\cE'$
whose objects are in both $\cE'$ and $\cE''$.
Notice that there is a chain of inclusions
\[ \cE \binc \cE'' \binc \cE''' \inc \cE'.
\]
We claim that each of these inclusions induces a weak equivalence on
nerves.  This is easy for $\cE''\inc \cE$ and $\cE'''\inc \cE'$, just
using functorial factorizations.

Define a functor $\theta\co \cE'' \ra \cE'''$ by sending an
object $B\cof E$ in $\cE''$ to the object $C \cof C\amalg_B E$.
To see that this lies in $\cE'''$ one can use Ken Brown's
lemma~\cite{hovey}
to show that pushing out a weak equivalence along a
cofibration yields another weak equivalence, provided all the domains
and codomains of the original maps are cofibrant.   It is simple to see
that $\theta$ gives a homotopy inverse for the inclusion $\cE'''\inc
\cE''$.

If $C\ra X$ lies in $\cE'$, a straightforward argument shows that $X$ is
weakly equivalent to $C\Wedge \Sigma^{n+1}M$ in the category $(C\ovcat R$--$\Alg
\ovcat C)$ (basically, one uses obstruction theory to directly
construct a zig-zag of weak equivalences).  So $\cE'$ is simply the
moduli space $\cM_{(C\ovcat R\text{--}\Alg \ovcat C)}(C\Wedge \Sigma^{n+1}M)$.

Recall from \fullref{prop.KI}
that there is a Quillen
equivalence
\[  \NonU_C \we (C\ovcat R\text{--}\Alg\ovcat C)
\]
in which the right adjoint $I$ sends $C\ra X \ra C$ to the fiber of
$X\ra C$.
Note that $I(C \Wedge \Sigma^{n+1}M) \simeq \Sigma^{n+1}M$ where
$\Sigma^{n+1}M$ is given the trivial structure of nonunital
$C$--algebra (in which the product is zero).
Thus, this Quillen equivalence implies that
$\cM_{(C\ovcat R\text{--}\Alg \ovcat C)}(C\Wedge
\Sigma^{n+1}M) \he \cM_{\NonU_C}(\Sigma^{n+1}M)$.
\end{proof}

We next reduce from $\NonU_C$ to $C\bimod$:

\begin{prop}
\label{pr:prop2b}
The functor $Z$ induces a weak equivalence
\[ \cM_{C\bimod}(\Sigma^{n+1}M) \ra \cM_{\NonU_C}(\Sigma^{n+1}M).
\]
\end{prop}

\begin{proof}
Since $Z$ preserves all weak equivalences, it induces a functor
between moduli categories in the obvious way.  Consider the
composite functor $\NonU_C \ra C\bimod$ given by
\[ X \mapsto P_{n+1}(Q(cX)), \]
where $c$ is any cofibrant-replacement functor in $\NonU_C$.  Applying 
\fullref{pr:Qanalysis} with $n$ replaced by $n+1$, we know 
that if $X$ is weakly equivalent to
$Z(\Sigma^{n+1}M)$ then $P_{n+1}(QcX)$ is a $C$--bimodule whose only
nonvanishing homotopy group lies in dimension $n+1$ and is isomorphic
to $M$.  So $P_{n+1}Qc$ induces a functor
$$F\co \cM_{NonU_C}(\Sigma^{n+1}M) \ra \cM_{C\bimod}(\Sigma^{n+1}M).$$

\fullref{pr:Qanalysis} implies that there is a natural
zig-zag of weak equivalences between the composite $F\circ Z$ and the
identity functor.  \fullref{pr:Qanalysis2} implies the
same for the composite $Z\circ F$.  It follows that the maps induced
by $Z$ and $F$ on the nerves of the categories are
homotopy inverses.
\end{proof}

We now need to analyze
$\cM_{C\bimod}(\Sigma^{n+1}M)$.  This is something which boils down to
an explicit computation.  For the next proposition, recall that to any
element $X$ of $\cM_{C\bimod}(\Sigma^{n+1}M)$ we may associate its
homotopy group $\pi_{n+1}X$ regarded as a $\pi_0(C)$--bimodule.  This
bimodule is isomorphic to $M$, and
in this way we obtain a functor 
$\pi_{n+1}\co \cM_{C\bimod}(\Sigma^{n+1}M) \ra \cM_{\pi_0(C)\bimod}(M)$.
The codomain is the category of $\pi_0(C)$--bimodules which are
isomorphic to $M$, with maps the isomorphisms; said
differently, it is the moduli space in the model category of
$\pi_0(C)$--bimodules where the weak equivalences are isomorphisms and
where every map is both a cofibration and a fibration.   We use this
functor in \fullref{7.8} and \fullref{thm:classify} to identify
the action of $\Aut(M)$.

\begin{prop}
\label{pr:prop2c}
Assume $C$ is cofibrant as an $R$--module (for example, $C$ is a
cofibrant $R$--algebra).  Then
the functor $\pi_{n+1}\co \cM_{C\bimod}(\Sigma^{n+1}M) \ra
\cM_{\pi_0(C)\bimod}(M)$
is a weak equivalence.  Consequently, one has
\[ \cM_{C\bimod}(\Sigma^{n+1}M) \he B\Aut(M)
\]
 where the automorphism
group is taken in the category of $\pi_0(C)$--bimodules.
\end{prop}

\begin{proof}
The Dwyer--Kan result \cite[2.3]{classifi}
%prop
identifies $\cM_{C\bimod}(\Sigma^{n+1}M)$ with the space
$B\hAut(\Sigma^{n+1}M)$,
where $\hAut$ denotes the simplicial monoid of homotopy automorphisms
in the model category $C\bimod$.  This is the subcomplex of the
homotopy function complex
$\underline{C\bimod}(\Sigma^{n+1}M,\Sigma^{n+1}M)$ consisting of all
path components which are invertible in $\pi_0$.  We'll now compute
this homotopy function complex.

The model category of $C$--bimodules is enriched, tensored, and
cotensored over symmetric spectra.  So for
any bimodules $N_1$ and $N_2$ there is a symmetric spectrum mapping
space $F(N_1, N_2)$.  The homotopy function complex is simply
the zero-th space $\Ev_0 F(cN_1,fN_2)$.  Also, since $C\bimod$ is a stable
model category one has $F(\Sigma cN_1,f\Sigma c N_2)\he F(cN_1,fN_2)$.
We need to use $F(c\Sigma^{n+1}M, f\Sigma^{n+1}M)\he F(cf M,cf M)$.

It is simple to compute that $\pi_i F(cf M,cf M)=0$ for $i>0$ and $\pi_0
F(cf M,cf M)=\Hom_{\pi_0(C)}(M,M)$, the group of endomorphisms of $M$ as a
$\pi_0(C)$--bimodule.  One way to do this is to recall that for any
two $C$--bimodules $N_1$ and $N_2$ there is a spectral sequence
\[ E_2^{p,q}=
\Ext^{p}_{\pi_*(C\Smash_R C^{op})}(\pi_*(N_1),\Sigma^{-q} \pi_*(N_2))
\Rightarrow \pi_{q-p} F(cN_1,fN_2).
\]
In the case $N_1=N_2=M$ the $E_2$--term completely vanishes in the
range $q\geq 0$ except for the single group $\smash{\Ext^0_{\pi_*(C\Smash_R
C^{op})}}(M,M)$ when $p=q=0$ (this uses that $C\Smash_R C$ is
connective, which in turn uses our cofibrancy assumption on $C$).
This group is the same as $\Hom_{\smash{\pi_0(C\Smash_R C^{op})}}(M,M)$.
Finally, we note that $\pi_0(C\Smash_R C^{op})\iso
\pi_0(C)\tens_{\pi_0 R} \pi_0(C)^{op}$.  This follows from the
spectral sequence
\[ \Tor_{p,q}^{\pi_*R}(\pi_*C,\pi_*C^{op})\Rightarrow \pi_{p+q}(C\Smash_R
C^{op})
\]
(again using our cofibrancy assumption on $C$),
together with the fact that $R$ and $C$ are connective.

Putting this all together, it readily follows that
$\hAut(\Sigma^{n+1}M)\he \Aut(M)$ (the latter regarded as a discrete
group).

Finally, we return to 
$\pi_{n+1}\co \cM_{C\bimod}(\Sigma^{n+1}M) \ra
\cM_{\pi_0(C)\bimod}(M)$.  Since the homotopy groups of both the
domain and codomain vanish except for $\pi_1$, it suffices to show the
functor induces an isomorphism on $\pi_1$.  Note that there are
obvious maps
\[ \Aut(M) \ra \pi_1 \cM_{C\bimod}(\Sigma^{n+1}M) \quad\text{and}\quad
\Aut(M) \ra \pi_1 \cM_{\pi_0(C)\bimod}(M).
\]
The first, for instance, sends an automorphism $\sigma$ to the loop
represented by the induced map of bimodules $\sigma
\co\Sigma^{n+1}M \ra \Sigma^{n+1}M$; the second is defined
similarly.  These maps obviously commute with the functor $\pi_{n+1}$.
But the map $\Aut(M) \ra \smash{\pi_1\cM_{\pi_0(C)\bimod}}(M)$ is readily seen
to be an isomorphism, and our analysis above of
$\hAut(\Sigma^{n+1}M)$ shows that the corresponding map
$\Aut(M) \ra \pi_1 \cM_{C\bimod}(\Sigma^{n+1}M)$ is also an isomorphism.
This finishes the proof.
\end{proof}

To any object $E\la B$ in $\cE_C(M,n)$ we may associate the abelian
group $\pi_{n+1}E$ which will be a $\pi_0(C)$--bimodule via the
isomorphism $\pi_0(B)\iso \pi_0(C)$ and the map $\pi_0(B)\ra
\pi_0(E)$.  So we have a functor
$\pi_{n+1}\co \cE_C(M,n) \ra
\cM_{\pi_0(C)\bimod}(M)$.

\begin{cor}
\label{7.8}
The functor $\pi_{n+1}\co \cE_C(M,n) \ra
\cM_{\pi_0(C)\bimod}(M)$ is a weak equivalence.  Consequently,
$\cE_C(M,n)\he B\Aut(M)$.
\end{cor}

\begin{proof}
Let $cC\we C$ be a cofibrant-replacement in the category of
$R$--algebras.  By \fullref{pr:changebase2}
the evident map $\cE_{cC}(M,n)\ra \cE_C(M,n)$ is a weak
equivalence.  We then have a zig-zag of weak equivalences
\[ \cE_C(M,n) \he \cE_{cC}(M,n) \he \cM_{NonU_{cC}}(\Sigma^{n+1}M) \he
\cM_{cC\bimod}(\Sigma^{n+1}M)
\]
provided by \fullref{lem:lem2b} and \fullref{pr:prop2b}.  
There is an obvious $\pi_{n+1}$ functor from each of these categories 
landing in $\cM_{\pi_0(C)\bimod}(M)$,
and the relevant triangles all commute.
Since $\pi_{n+1}\co \cM_{cC\bimod}(\Sigma^{n+1}M) \ra
\smash{\cM_{\pi_0(C)\bimod}}(M)$ is a weak equivalence by
\fullref{pr:prop2c}, we deduce that the same is true for the
$\pi_{n+1}$ functor with domain $\cE_C(M,n)$.
\end{proof}

We have finally completed our main proof:

\begin{proof}[Proof of \fullref{pr:prop2}]
Part (a) follows directly from \fullref{lem:lem2a} and the remarks
following its proof.
Part (b) is a consequence of the preceding corollary.
\end{proof}

%%%%%%%%%%%%%%%%%%%%%%%%%%%%%%%%%%%%%%%%%%%%%%%%%%%%%%%

\section{The main result}

Finally, we can pull everything together and prove the main theorem:

\begin{thm}
\label{thm:classify}
Fix $n\geq 1$.
Let $R$ be a connective ring spectrum, and let $C$ be a connective
$R$--algebra such that $P_{n-1}C\he C$.  Let $M$ be a $\pi_0(C)$--bimodule.
There is a homotopy fiber sequence of spaces
\[ \underline{R\text{--}\Alg_{\smash{/C}}}(C,C\Wedge\Sigma^{n+1}M) \ra |\cM(C{+}(M,n))| \ra
B\Aut(M).
\]
Consequently, one has a bijection
\[ [\Ho(R\text{--}\Alg_{/C})(C,C\Wedge\Sigma^{n+1}M)]/\Aut(M) \iso \pi_0
\cM(C{+}(M,n))
\]
where $\Aut(M)$ acts on the second factor of $C \Wedge\Sigma^{n+1}M$. 
\end{thm}

\begin{proof}
First one uses that $\cM(C{+}(M,n))\he \cK_C(M,n)$, from
\fullref{pr:prop1}.  Then one uses the homotopy fiber sequence
\[ \cM_{R\text{--}\Alg_{/C}}(C,C\Wedge\Sigma^{n+1}M) \ra \cK_C(M,n) \ra
\cE_C(M,n)
\]
established in \fullref{lem:lem2a} and our identification
$\cE_C(M,n)\he B\Aut(M)$.  This proves the first claim of the theorem.

To prove the second statement, we look at the long exact homotopy
sequence for the above fiber sequence.  Since $\cE_C(M,n)$ is
connected, this identifies $\smash{\pi_0\cK_C}(M,n)$ with a quotient of
$\pi_0[\cM_{R\text{--}\Alg_{/C}}(C,C\Wedge\Sigma^{n+1}M)]$ by an action of
$\pi_1\cE_C(M,n)\iso \Aut(M)$.  We must identify the action.

For brevity, write $S=\pi_0[\cM_{R\text{--}\Alg_{/C}}(C,C\Wedge\Sigma^{n+1}M)]$.
Every equivalence class $s\in S$ can be represented by
a diagram
\[ \xymatrix{C & A \ar[l]_{\hbox{\footnotesize $\sim$}} \ar[r]^-{g} & f(C\Wedge \Sigma^{n+1}M)
& C\Wedge\Sigma^{n+1}M. \ar[l]}
\]
Let $\sigma\in \Aut(M)$.  Under the identification
$\pi_1\cE_C(M,n)\iso \Aut(M)$, $\sigma$ corresponds to the self-map of
the object $[C\inc C\Wedge\Sigma^{n+1}M]$  which is the identity on
$C$ and induced by $\sigma$ on $M$.  We can just as well represent
$\sigma$ as a self-map of $[C\inc f(C\Wedge\Sigma^{n+1}M)]$.

To determine the action of $\sigma$ on $s$ we do the usual thing: we
lift the loop represented by $\sigma$ to a path in $\cK_C(M,n)$
beginning at $s$, and we take the terminal point of that path.
Our path is the map
\[ \xymatrix{ A \ar[r]^-g\ar@{=}[d] & f(C\Wedge
\Sigma^{n+1}M)\ar[d]_\sigma
& C\ar@{ >->}[l] \ar@{=}[d]
\\
A \ar[r]^-{\sigma g} & f(C\Wedge \Sigma^{n+1}M) & C\ar@{ >->}[l]
}
\]
This identifies the action of $\Aut(M)$ on
$S=\Ho(R$--$\Alg_{/C})(C,C\Wedge \Sigma^{n+1}M)$ with the action coming
from the second factor of $C\Wedge\Sigma^{n+1}M$.
\end{proof}

\begin{remark}
The isomorphism from $[\Ho(R$--$\Alg_{/C})(C,C\Wedge\Sigma^{n+1}M)]/\Aut(M)$
to $\pi_0 \cM(C{+}(M,n))$ produced in the theorem is precisely the
pullback map $\widetilde{\PB}$ defined in
\fullref{se:exten-pullback}.
This follows at once by looking at the various maps we used in our
identifications (particularly the one of \fullref{pr:prop1}).
\end{remark}

\begin{remark}
\fullref{thm:main} and \fullref{cor:main} from the
introduction are the special case of the above theorem where $R=S$.
\end{remark}

\begin{remark}
\label{re:bdg}
In~\cite{BDG} the moduli problem for Postnikov extensions of spaces is
set up a bit differently.  Proposition 3.7 in \cite{BDG} considers all spaces
$Y\simeq P_n Y$ for which there exists some chain of weak equivalences
$P_{n-1} Y \simeq C$, instead of considering maps $Y \to C$ with a
fixed space $C$ as the $(n-1)$st Postnikov section.  (Note that with the
definition of \cite{BDG} one must be careful about the isomorphism
$\pi_{n+1}Y \iso M$, as this must be an isomorphism of bimodules and
to make sense of this one needs a fixed
isomorphism of rings $\pi_0Y\iso \pi_0C$). 

If one adopts the choices of \cite{BDG}, one obtains a homotopy fiber
sequence analogous to \fullref{thm:classify} where the base space
is a product $\cM(C)\times B\Aut(M)$ and $\cM(C)$ is the moduli space
of $R$--algebras weakly equivalent to $C$.  See \cite[3.10]{BDG}.  It
is probably possible to directly relate our \fullref{thm:classify} to
this formulation, although we have not pursued this.  
\end{remark}

\subsection[Alternative formulations and the connection with THH]{Alternative formulations and the connection with $\THH$}
\label{sec:thh}
Let $R$, $C$, $M$, and $n$ be as in the statement of
\fullref{thm:classify}.  In a typical application, one is
interested in the possible extensions of $C$ by $M$ and thus is led to
try to compute $\Ho(R$--$\Alg_{/C})(C,C\Wedge \Sigma^{n+1}M)$.  There are two
well-known ways to simplify this, which we briefly record here.

In some applications one is interested in the following reduction.
We don't use it in the present paper, but it makes sense to record it here. 

\begin{prop}\label{prop-PC}
Assume $C$ is a cofibrant $R$--algebra and fix a zero-th
Postnikov section $p_0\mc C \to P_0 C$.  One has a weak equivalence of
homotopy function complexes
\[
\underline{R\text{--}\Alg_{\smash{/C}}}(C,C\Wedge\Sigma^{n+1}M) \he
\underline{R\text{--}\Alg_{\smash{/P_0C}}}(C,P_0C \Wedge\Sigma^{n+1}M).
\] 
\end{prop}

\begin{proof}
First note that one has a Quillen pair $L\co R$--$\Alg_{/C} \adjoint
R$--$\Alg_{/P_0C}\co T$ where $L$ is composition with $p_0$ and $T$ is
base-change along $p_0$.  This induces a weak equivalence of homotopy
function complexes
\[ \underline{R\text{--}\Alg_{\smash{/P_0C}}}(L(C),P_0C\Wedge\Sigma^{n+1}M) \he
\underline{R\text{--}\Alg_{\smash{/C}}}(C,\underline{T}(P_0C\Wedge\Sigma^{n+1}M)) 
\]
where $\underline{T}$ denotes the derived functor of $T$ (and where we
have used that $C$ is cofibrant).  Thus, it suffices for us to show
that
$\underline{T}(P_0C\Wedge\Sigma^{n+1}M) \he C\Wedge \Sigma^{n+1}M$ in
$R$--$\smash{\Alg_{/C}}$.  

Factor the projection $P_0C \Wedge \Sigma^{n+1}M \to P_0C$
into a trivial cofibration followed by a fibration
$P_0C \Wedge\Sigma^{n+1}M \smash{\stackrel{i~}{\vphantom{-}\smash{\rightarrow}}} \fP \smash{\stackrel{h~}{\vphantom{-}\smash{\rightarrow}}} P_0 C$.
Define $\fC$ as the pullback in the following square
\[ \xymatrix{ 
\fC \ar[r]^{p_0'}\ar@{>>}[d]^{h'} & \fP \ar@{>>}[d]^{h} \\ 
C \ar[r]^{p_0} & P_0C }
\]
and note that $\fC$ is a model for $\underline{T}(P_0C\Wedge \Sigma^{n+1}M)$.  

Observe that the map $C \Wedge \Sigma^{n+1}M \to P_0C \Wedge \Sigma^{n+1}M
\to \fP$ factors through $\fC$.
To see that $C \Wedge \Sigma^{n+1}M \to \fC$
is a weak equivalence note that in the category of spectra the homotopy fiber 
of both $\fC \to C$ and $C \Wedge \Sigma^{n+1}M \to C$ is $\Sigma^{n+1}M$.  
\end{proof}

The second well-known simplification is given by  
a connection with topological Hoch\-schild cohomology, $\THH^*$.  
This is via three Quillen pairs with the left adjoints on top:
\[\xymatrix@=10mm{
R\text{--}\Alg_{/C}  \ \ %\quad
\ar@<.4ex>^-{F}[r] &
%\quad  
\ \  (C\ovcat R\text{--}\Alg\ovcat C) \ \    %\quad
\ar@<.4ex>^-{G}[l]
\ar@<-.4ex>_-{I}[r] &
%\quad 
 \ \  \NonU_C \ \ %\quad
\ar@<-.4ex>_-{K}[l]
\ar@<.4ex>^-{Q}[r] &
\ \ % \quad 
C\bimod %\quad 
\ar@<.4ex>^-{Z}[l]
}
\]
%
%\[ F\co R\text{--}\Alg_{/C} \adjoint (C\ovcat R\text{--}\Alg\ovcat C)\co G, \qquad
%K\co \NonU_C \adjoint
%(C\ovcat R\text{--}\Alg\ovcat C)\co I, 
%\]
%\[ \text{and} \qquad Q\co \NonU_C \adjoint C\bimod \co Z.
%\]
%
The second two of these were defined in \fullref{sec:nonu}.  In
the first, we have $F(X)=C\amalg_R X$ and $G$ is the forgetful
functor.  Recall that the $(K,I)$ pair is a
Quillen equivalence by \fullref{prop.KI}.

Following \cite{B} (but using slightly different notation), for $X\in
R$--$\smash{\Alg_{/C}}$ one defines
\[ \Omega_{R\ra C}(X)=(
\underline{Q} \circ \underline{I} \circ \underline{F})(X). \]
Recall that the underlines denote derived functors.
One obtains a chain of weak equivalences
\begin{align*}
\underline{R\text{--}\Alg_{\smash{/C}}}(C,C\Wedge\Sigma^{n+1}M) & \he
\underline{(C\ovcat R\text{--}\Alg\ovcat
C)}\, \bigl (\underline{F}(C),C\Wedge\Sigma^{n+1}M \bigr )\\
& \he
\underline{\NonU_C}\bigl (\underline{I}(\underline{F}C), \Sigma^{n+1}M
\bigr ) \\
& \he
\underline{C\bimod}(\Omega_{R\ra C}(C),\Sigma^{n+1}M).
\end{align*}
In this chain we have used three facts, namely 
\[
 \underline{G}(C\Wedge \Sigma^{n+1}M)\he C\Wedge \Sigma^{n+1}M, %\mbox{,  \  }
\qquad
\underline{I}(C\Wedge \Sigma^{n+1}M)\he \Sigma^{n+1}M, 
\]
\[
\underline{Z}(\Sigma^{n+1}M)\he \Sigma^{n+1}M. \leqno{\hbox{and}}
\]
The first and last are trivial.  For the second, use the fact that the
evident map from $\Sigma^{n+1}M$ to the homotopy fiber of
$C\Wedge\Sigma^{n+1}M \ra C$ is a weak equivalence.

Assume that $C$ is cofibrant as an $R$--module, and recall that
$C\bimod$ denotes the category of $(C\Smash_R C^{op})$--modules.
One defines the group of derivations
\begin{align*}
 \Der^{n+1}_R(C,M) &=\ho({C\bimod})(\Omega_{R\ra
C}(C),\Sigma^{n+1}M)\\
&=\pi_{-n-1}\Bigl[ \underline{C\bimod}(\Omega_{R\ra
C}(C), M) \Bigr ].   %%%%%%%%%%%  ,
\end{align*}
Note that this is just $\pi_0$ of $\underline{R\text{--}\Alg_{\smash{/C}}}(C,C\Wedge
\Sigma^{n+1}M)$  by the weak equivalence of mapping spaces given
above.

One also defines 
\[
\THH^k_R(C,M)=\ho({C\bimod})(C,\Sigma^{k}M)=\pi_{-k}\Bigl [
\underline{C\bimod}(C,M)\Bigr ].
\]  
To connect these groups, one identifies $\Omega_{R\ra C}(C)$ with the
homotopy fiber of the multiplication map
\[ C \Smash_R C^{op} \ra C\]
(this uses that
$C$ is cofibrant as an $R$--module, otherwise a
cofibrant-replacement is necessary before forming the smash product).
This identification is nontrivial, but a proof has been shown to us by
Mike Mandell \cite{M}.

\begin{remark}
In \cite{L} the bimodule $\Omega_{R\ra C}(C)$ is {\it defined\/} to be
the homotopy fiber of the above multiplication map.  The hard work is
then to prove that the mapping spaces
$\underline{C\bimod}(\Omega_{R\ra C}(C),\Sigma^{n+1}M)$ and
$\underline{R\text{--}\Alg_{\smash{/C}}}(C,C\Wedge\Sigma^{n+1}M)$ are weakly
equivalent.  The proof of this fact in \cite{L} contains gaps.
\end{remark}

Applying $\underline{C\bimod}(\blank,\Sigma^{n+1}M)$ to the homotopy
fiber sequence $\Omega_{R\ra C}(C) \ra C\Smash_R C \ra C$ induces a
homotopy fiber sequence of mapping spaces.  Consider the associated
long exact homotopy sequence.  One has $\pi_i\big(
\underline{C\bimod}(C\Smash_R C,\Sigma^{n+1}M) \big) \iso
\pi_i(\Sigma^{n+1}M)$; so as long as $n\geq 0$, this group vanishes
for $i\leq 0$.  The long exact sequence then shows that for $n\geq 0$,
\[ \Der^{n+1}_R(C,M)\iso \THH^{n+2}_R(C,M).
\]
Putting everything together, we have proven the following:

\begin{prop}\label{prop-8.8}
Let $R$, $C$, $M$, and $n$ be as in \fullref{thm:classify}.
Assuming that $C$ is cofibrant as an $R$--module, one has 
a bijection
\[ \THH^{n+2}_R(C,M)/\Aut(M) \iso \pi_0\cM(C{+}(M,n)).\]
\end{prop}

\begin{remark}
In the case where $C$ is not cofibrant as an $R$--module, one can
repeat the above discussion by replacing $C\bimod$ with $QC\bimod$,
where $QC\ra C$ is a cofibrant-replacement for $C$ in $R$--$\Alg$.  The
correct definition of $\THH^i_R(C,M)$ should really be
$\ho(QC\bimod)(QC,\Sigma^i M)$, and analogously for $\Der$. 
\end{remark}

%%%%%%%%%%%%%%%%%%%%%%%%%%%%%%%%%%%%%%%%%%%%%%%%%

\appendix

\section{Proof of the Blakers--Massey theorem for ring spectra}
\setobjecttype{App}
\label{se:pushout-proof}

In this section we will give the proof of
\fullref{thm:smallBM}.  To ease notation, refer to a map $f\mc
X \to Y$ as an $(n-1)$--equivalence if $\pi_i(Y,X) = 0$ for $i < n$.
Recall our conventions that we replace $f$ by a cofibration before
considering relative homotopy, and that homotopy groups always refer
to the derived homotopy groups ($\pi_*$ applied to a fibrant
replacement).

We actually prove the following result.  For this statement let
$\amalg$ denote the coproduct of $R$--algebras and $\cup$ denote the
coproduct of $R$--modules.  For example, in this notation the homotopy
pushout appearing in \fullref{thm:smallBM} is $P = C
\amalg^h_{A} B$.

\begin{thm}
\label{thm:actualBM} 
Let $R$ be a connective, commutative ring spectrum and let $m,n \geq 1$.
Suppose given two maps of $R$--algebras: 
$A\to B$ an $(n-1)$--equivalence and $A \to C$ an $(m-1)$--equivalence,
with $A$ connective.
Then the map from the homotopy pushout of
$R$--modules to the homotopy pushout of $R$--algebras, $C\cup^h_A B \to
C\coprd^h_A B$, is an $(m+n-1)$--equivalence.
\end{thm}

We prove this proposition at the end of this section.  To get an idea
of why it should be true, it's useful to think about the analogous
result for dgas.  There, by replacing $B$ up to quasi-isomorphism one
can assume that $A\ra B$ is a monomorphism which is the identity in
degrees strictly smaller than $n$.  Similarly, one can assume $A\ra C$
is a monomorphism which is the identity in degrees strictly smaller
than $m$.  The coproduct $B\coprd_A C$ is constructed from formal
words in the elements of $B$ and $C$, and inspection shows immediately
that $B\cup_A C \ra B\coprd_A C$ is an isomorphism in degrees less
than $m+n$.  Our proof for the ring spectra case will have a similar
flavor, but it must contend with the fact that ring spectra do not
have ``elements''.

\begin{remark}
Throughout this section we need to deal with colimits both of
$R$--modules and $R$--algebras.  There are two basic facts which will be
used repeatedly.  First, colimits in the category of $R$--modules are
the same as colimits in the category of $S$--modules (ie in the
category of spectra).  Said better, the underlying spectrum of a
colimit of $R$--modules is the same as the colimit of the underlying
spectra.

The second fact we will need is that for a sequential colimit of
$R$--algebras, the colimit in the category of $R$--algebras is the same
as the colimit in the category of $R$--modules.  Both of these facts
work quite generally in the context of arbitrary symmetric monoidal
categories.
\end{remark}

Before tackling the proof of \fullref{thm:actualBM}, we show
how the result implies the Blakers--Massey theorem:

\begin{proof}[Proof of \fullref{thm:smallBM}]
We can assume that $A, B$ and $C$ are cofibrant and fibrant as
$R$--algebras (and hence fibrant as spectra).  We can also assume $A
\to B$ and $A \to C$ are cofibrations of $R$--algebras and hence also
cofibrations of $R$--modules by~\cite[4.1(3)]{SS1}.  Then the homotopy
pushout of $R$--modules is the pushout $C \cup_A B$ and the homotopy
pushout of $R$--algebras is the pushout $P= C \coprd_A B$.

By \fullref{thm:actualBM}, we know that $C\cup_A B \to
C\coprd_A B$ is an $(m+n-1)$--equivalence.  It follows that
$\pi_i(C\cup_A B, C) \to \pi_i(P, C)$ is an isomorphism for $i <
m+n-1$ and a surjection for $i = m+n-1$.  But $\pi_i(C\cup_A B,C)\iso
\pi_i((C\cup_A B)/C)\iso \pi_i(B/A)\iso \pi_i(B,A)$.  So $\pi_i(B,A)
\to \pi_i(P,C)$ is an isomorphism for $i < m+n-1$ and a surjection for
$i = m+n-1$.
\end{proof} 

The proof of \fullref{thm:actualBM} will require several
lemmas.  Although the next lemma is stated for an arbitrary
$(n-1)$--connected cofibration $K\cof L$ of pointed simplicial sets,
the main application is when $K\ra L$ is $\bdd{n}\ra\del{n}$ or a
coproduct of such maps.  Let $\pi^s_i(L,K)$ denote
$\pi_i(\Sigma^{\infty}L, \Sigma^{\infty}K)$.  Recall that for a
pointed simplicial set $K$, $T_R(K)$ is our shorthand for the tensor
algebra $T_R(R\Smash \Sigma^\infty K)$.

\begin{lemma}
\label{lem:smallBM}
Let $n\geq 2$ and 
let $R$ be a connective, commutative ring spectrum.
Let $K\cof L$ be
a cofibration of pointed simplicial sets such that $\pi^s_i(L,K)=0$ for
$i< n$.
Suppose $X$ is a
cofibrant, connective $R$--algebra and
\[ \xymatrix{
T_R(K) \ar[r]\ar[d] & X \ar[d] \\
T_R(L) \ar[r] & Y
}
\]
is a pushout square of $R$--algebras.
Then $X\ra Y$ is an $(n-1)$--equivalence.
\end{lemma}

\begin{proof}
In \cite[Proof of Lemma 6.2]{SS1} the pushout $Y$ is described as a
certain directed colimit of pushouts in the category of $R$--modules.
If we let $P_0=X$, then there is a sequence of cofibrations of $R$--modules 
\[ P_0 \cof P_1 \cof P_2 \cof \cdots
\]
whose colimit is the underlying $R$--module of $Y$ and where $P_{r}$
is obtained from $P_{r-1}$ by a pushout diagram of $R$--modules 
\[ \xymatrix{
W_r \Smash (X\Smash_R X\Smash_R \cdots \Smash_R X) \ar[r]\ar[d] & P_{r-1}
\ar[d] \\
W'_r \Smash (X\Smash_R X\Smash_R \cdots \Smash_R X) \ar[r] & P_{r}.
}
\]
Here there are $r+1$ copies of $X$ in the smash product, and $W_r \ra
W'_r$ is the $r$--fold box product of $K \ra L$.  So $W'_r/W_r \iso
(L/K)^{\Smash r}$, which is $(rn-1)$--connected.

Note that $X\Smash_R X\Smash_R\cdots \Smash_R X$ is connective, since
both $X$ and $R$ are connective and $X$ is cofibrant.  This follows 
from the fact that $X\Smash_R X$
is the realization of the simplicial $R$--module $[n]\mapsto X\Smash
R^{\Smash n}\Smash X$, for instance.  So we find that
the $R$--module 
\[P_r/P_{r-1}\iso (L/K)^{\Smash r}\Smash (X\Smash_R\Smash\cdots
\Smash_R X)
\]
has no homotopy groups below dimension $rn$.
Hence, $P_{r-1}\ra P_r$ is an $(rn-1)$--equivalence.  
It follows immediately that $X \to Y$ is an $(n-1)$--equivalence.
\end{proof}

Consider the following subset of the generating cofibrations
of $R$--algebras: $I'_n =\{ T_R(\bdd{l})\ra T_R(\del{l}) \,|\, l \geq n\}$.
Let $I_n = I'_n \cup J$ where $J$ is the set of generating trivial
cofibrations of $R$--algebras. 
Note that each of these maps is an $(n-1)$--equivalence.
The next lemma converts any $(n-1)$--equivalence of $R$--algebras
into a weakly equivalent map built from $I_{n}$ by colimits and pushouts. 
Let $I_{n}\Cell$ be the collection of maps which are (possibly infinite)
compositions of pushouts of maps in $I_{n}$~\cite[2.1.9]{hovey}. 

\begin{lemma}
\label{lem-nc}
If $f\mc A \to B$ is an $(n-1)$--equivalence 
between fibrant $R$--algebras, then there is a factorization
$A \smash{\stackrel{i~}{\vphantom{-}\smash{\rightarrow}}} B' \smash{\stackrel{p~}{\vphantom{-}\smash{\rightarrow}}} B$ with $f = pi$, $i \in I_{n}\Cell$
and $p$ a trivial fibration.  
\end{lemma}

\begin{proof}
The small object argument produces a factorization of $f$  as
$f=pi$ such that $i\mc A \to B'$ is in $I_{n}\Cell$ and $p\mc B' \to B$ has 
the right lifting property with respect to $I_{n}$.  
Since the maps in $J$ are trivial cofibrations, any map in $J\Cell$ is
a weak equivalence.  Thus, by \fullref{lem:smallBM},  $i$ is a 
(possibly infinite) composition of $(n-1)$--equivalences. 
So $i$ is
an $(n-1)$--equivalence.  It follows that $p$ is an 
$(n-1)$--equivalence as well;  hence $\pi_i(B,B')=0$ for $i < n$. 
The map $p$ is a fibration since it has the right lifting property
with respect to $J$; hence, $B$ and $B'$ are both fibrant.
Using the definition in \fullref{se:relative} of relative homotopy 
groups, it follows 
that $\pi_i(B,B')=0$ for $i \geq n$ as well since $p$ has the right lifting 
property with respect to $I'_{n}.$  Thus, $p$ is a trivial fibration. 
\end{proof}

\begin{lemma}
\label{lem-tinyBM}
Assume $R, K$ and $L$ satisfy the hypotheses of \fullref{lem:smallBM}.
Suppose $X$ and $C$  are cofibrant, connective $R$--algebras, 
$X \to C$ is an $(m-1)$--equivalence and the following two squares
\[ \xymatrix{
T_R(K) \ar[r]\ar[d] & X \ar[d]\ar[r] & C \ar[d] \\
T_R(L) \ar[r] & Y \ar[r] & Q
}
\]
are both pushout squares of $R$--algebras.
Then the map from the coproduct in $R$--modules to the coproduct
of $R$--algebras, $C \cup_X Y \to C\coprd_X Y = Q$, is an $(m+n-1)$ 
equivalence.  
\end{lemma}

Note here $Y= X \coprd_{T_R(K)} T_R(L)$ and $Q= C \coprd_X Y 
= C \coprd_{T_R(K)} T_R(L)$.  

\begin{proof}
As described in the proof of \fullref{lem:smallBM},
$Y$ is the colimit of a sequence of cofibrations of $R$--modules 
\[ X =P_0 \cof P_1 \cof P_2 \cof \cdots
\]
with $P_r/P_{r-1} \iso (L/K)^{\sm r} \sm (X \sm_R \cdots \sm_R X)$.
Similarly, one can
produce $Q$ as a colimit of a sequence of cofibrations of $R$--modules 
\[ C =Q_0 \cof Q_1 \cof Q_2 \cof \cdots
\]
with $Q_r/Q_{r-1} \iso (L/K)^{\sm r} \sm (C \sm_R \cdots \sm_R C)$.
The map $X \to C$ induces a map of sequences $P_i \to Q_i$.

Next we build a sequence $\smash{P'_i}$ between $P_i$ and $Q_i$. 
Set $\smash{P'_0}=C$ and define $\smash{P'_i}$ as the pushout of $R$--modules 
$\smash{P'_{i-1}} \cup_{P_{i-1}} P_i$;
so each square in the following diagram is a pushout of $R$--modules.
\[ \xymatrix{
X \ar@{=}[r]\ar[d] & 
P_0 \ar[r]\ar[d] & P_1 \ar[d]\ar[r] & P_2 \ar[d]\ar[r] & \cdots \\
C \ar@{=}[r] & P'_0 \ar[r] & P'_1 \ar[r] & P'_2 \ar[r] & \cdots 
}
\]
Since $\colim P_i = Y$, it follows formally that  
$\colim \smash{P'_i} = C \cup_X Y$.  Also note that the maps $P_i \to Q_i$
induce maps $f_i\mc \smash{P'_i} \to Q_i$.  

We next show inductively that each $f_i\mc \smash{P'_i} \to Q_i$ is an 
$(m+n-1)$--equivalence.  
First, $f_0 \mc C \to C$ is the identity, and hence an $(m+n-1)$
equivalence.  By construction one has $\smash{P'_r/P'_{r-1}} \iso P_r/P_{r-1} \iso
(L/K)^{\sm r} \sm (X \sm_R \cdots \sm_R X)$ and 
$Q_r/Q_{r-1} \iso (L/K)^{\sm r} \sm (C \sm_R \cdots \sm_R C)$.
Note that $X\sm_R \cdots \sm_R X \ra C\sm_R \cdots \sm_R C$ is the
composite
\[ X^{\sm_R(r+1)} \ra C\sm_R (X^{\sm_R(r)}) \ra C\sm_R C \sm_R
(X^{\sm_R(r-1)}) \ra \cdots \ra C^{\sm_R(r+1)}.
\]
Each map in this sequence comes from smashing $X\ra C$ over $R$ with a
connected (and cofibrant) $R$--module, and so each map is an
$(m-1)$--equivalence.  The spectrum $(L/K)^{\sm r}$ is
$(rn-1)$--connected, so
\[ (L/K)^{\sm r} \sm (X^{\sm_R(r+1)}) \ra
(L/K)^{\sm r} \sm (C^{\sm_R(r+1)}) 
\]
is an $(rn+m-1)$--equivalence.  Hence $P'_r/P'_{r-1} \ra
Q_r/Q_{r-1}$ is, in particular, an $(m+n-1)$--equivalence.

Since $\smash{P'_{r-1}} \to Q_{r-1}$ is an $(m+n-1)$--equivalence by induction,
it now follows that $\smash{P'_r} \to Q_r$ is also an $(m+n-1)$--equivalence.
So this holds for all $r$.
Since homotopy groups commute with colimits of $R$--modules,
this implies that $\colim P'_r \to \colim Q_r $ is an $(m+n-1)$--equivalence.
\end{proof}

Finally we can complete the main result of this section:

\begin{proof}[Proof of \fullref{thm:actualBM}]
We can assume that $A$, $B$, and $C$ are cofibrant and fibrant as
$R$--algebras (hence fibrant as spectra).  By \fullref{lem-nc}, 
we can assume $A\ra B$ is in $I_{n}\Cell$.  We are now in a situation
where the homotopy pushouts are weakly equivalent to the pushouts, so
our task is to show that $C\cup_A B \ra C\coprd_A B$ is an
$(m+n-1)$--equivalence. 

Since $A \ra B$ is in $I_{n}\Cell$ we can assume $B=\colim B_i$ with
$B_0 = A$ and each $B_{i-1} \to B_{i}$ a pushout of a coproduct of
maps in $I_{n}= I'_{n} \cup J$.  
Any coproduct of maps in $I_n'$ has the form $T_R(K)\ra T_R(L)$ for some
cofibration of pointed simplicial sets $K\ra L$ which is a stable
$(n-1)$--equivalence.  And any coproduct of maps in $J$ is a trivial
cofibration.  So any coproduct of maps in $I_n$ has the form
\[ T_R(K)\amalg W' \ra T_R(L)\amalg Z' \]
where $W'\ra Z'$ is a trivial cofibration.  Forming $B_i$ from
$B_{i-1}$ can therefore be done in two stages, by first pushing out
along the map $T_R(K)\ra T_R(L)$ and then pushing out along the map
$W'\ra Z'$.  By now redefining the $B_i$'s, we can assume that each
map $B_{i-1}\ra B_i$ in our colimit is obtained by pushing out either
along a trivial cofibration or along a map $T_R(K)\ra T_R(L)$ as
above.

We will show inductively that each map $f_i \mc C \cup_A B_i \to C
\coprd_A B_i$ is an $(m + n -1)$--equivalence and then conclude that
the colimit $C\cup_A B \to C\coprd_A B$ is also an
$(m+n-1)$--equivalence.  (Here we are using that a sequential colimit
of $R$--algebras is the same as the sequential colimit of underlying
$R$--modules.)

Note that $f_0\mc C \to C$ is an isomorphism since $B_0 = A$.  Assume
$f_{i-1}$ is a $(m+n-1)$--equivalence.  If $B_{i-1} \to B_i$ is a
pushout of a trivial cofibration, then it is a trivial cofibration of
$R$--algebras and hence also a trivial cofibration of $R$--modules
\cite[4.1(3)]{SS1}.  This uses the fact that $B_{i-1}$ is cofibrant
as an $R$--algebra and hence also as an $R$--module.  It follows that
both $C\cup_A B_{i-1} \to C\cup_A B_i$ and $C\coprd_A B_{i-1} \to
C\coprd_A B_i$ are weak equivalences.  Thus, if
$f_{i-1}$ is an $(m+n-1)$--equivalence then so is $f_i$.

For the remaining case 
we have $B_i = B_{i-1} \coprd_{T_R(K)} T_R(L)$ with $K\to L$ as in 
\fullref{lem:smallBM}.  We show that 
$f_i\mc C\cup_A B_i \to C \coprd_A B_i$ is the 
composition of two $(m+n-1)$--equivalences.  For the 
first piece, consider the diagram of pushout squares in $R$--modules: 
\[ \xymatrix{
B_{i-1} \ar[r]\ar@{ >->}[d] & C\cup_A B_{i-1} \ar@{ >->}[d]\ar[r]^{f_{i-1}} & C\coprd_A B_{i-1} 
\ar@{ >->}[d] \\
B_i \ar[r] & C\cup_A B_i \ar[r]^{f'_i \ \ \ \ \ \ \ } & (C\coprd_A B_{i-1})\cup_{B_{i-1}}B_i 
}
\]
The map $f'_i$ is an $(m+n-1)$--equivalence since it is the pushout of the 
$(m+n-1)$--equivalence $f_{i-1}$.  Next consider the pushout of $R$--algebras
\[\xymatrix{
B_{i-1} \ar[r]\ar[d] & C\coprd_A B_{i-1} \ar[d] \\
B_i \ar[r] & C\coprd_A B_i 
}
\]
We claim that the top map here is an $(m-1)$--equivalence.  Assuming
this and recalling that 
$B_i = B_{i-1} \coprd_{T_R(K)} T_R(L)$, then 
by \fullref{lem-tinyBM} that $\smash{f''_i}\mc (C\coprd_A B_{i-1}) 
\cup_{B_{i-1}} B_i \to C \coprd_A B_i$ is an $(m+n-1)$--equivalence.  Since
$f_i = \smash{f''_i f'_i}$ this would finish the induction step.
To establish the claim about the top map above, compare the two horizontal
cofibration sequences below. 
\[\xymatrix{
B_{i-1} \ar[r]\ar[d]^{\iso} & C\cup_A B_{i-1} \ar[d]^{f_{i-1}}\ar[r] & (C\cup_A B_{i-1})/B_{i-1} = C/A\ar[d] \\
B_{i-1} \ar[r]& C\coprd_A B_{i-1} \ar[r] & (C\coprd_A B_{i-1})/ B_{i-1} 
}
\]
\newpage
Since the left vertical map is an isomorphism and the middle map
is an $(m+n-1)$--equivalence the right map is also an $(m+n-1)$--equivalence.
Since $C/A$ is $(m-1)$--connected, so is $(C\coprd_A B_{i-1})/B_{i-1}$. 
It follows that $B_{i-1} \to C\coprd_A B_{i-1}$
is an $(m-1)$--equivalence.
\end{proof}

%%%%%%%%%%%%%%%%%%%%%%%%%%%%%%%%%%%%%%%%%%%%%%%%%%%%%%

\bibliographystyle{gtart}
\bibliography{link}

\end{document}